\documentclass{article}
\usepackage[utf8]{inputenc}

\usepackage{geometry}
\usepackage{amsmath}
\usepackage{amsfonts}
\usepackage{float}
\usepackage{stmaryrd}
\usepackage{graphicx}
\usepackage{caption}
\usepackage{subcaption}
\usepackage{authblk}
\usepackage[all]{xy}
\usepackage{tikz-cd}
\usepackage{hyperref}
\usepackage{cleveref}
\usepackage{multirow}
\usepackage{amsthm}
\usepackage{grffile}

\newcommand{\rhoh}{\rho_h}
\newcommand{\sh}{s_h}
\newcommand{\Bh}{\BB_h}
\newcommand{\Ah}{\AB_h}
\newcommand{\uh}{\uu_h}
\newcommand{\vh}{\vv_h}
\newcommand{\mh}{\mm_h}

\newcommand{\uu}{\boldsymbol{u}}
\newcommand{\mm}{\boldsymbol{m}}
\newcommand{\vv}{\boldsymbol{v}}

\newcommand{\AB}{\boldsymbol{A}}
\newcommand{\BB}{\boldsymbol{B}}

\newtheorem{theorem}{Theorem}
\newtheorem{proposition}[theorem]{Proposition}%

\newtheorem{problem}{Problem}%

\newtheorem{remark}{Remark}%

\newtheorem{definition}{Definition}

\crefname{problem}{problem}{problems}
\crefname{assumption}{assumption}{assumptions}
\title{Variational discretizations of viscous and resistive magnetohydrodynamics using structure-preserving finite elements}
\author{Valentin Carlier}
\date{\today}

\begin{document}

\maketitle

\begin{abstract}
We propose a novel structure preserving discretization for viscous and resistive magnetohydrodynamics. 
We follow the recent line of work on discrete least action principle for fluid and plasma equation, incorporating the recent advances to model dissipative phenomena through a generalized Lagrange-d'Alembert constrained variational principle.
We prove that our semi-discrete scheme is equivalent to a metriplectic system and use this property to propose a Poisson spliting time integration.
The resulting approximation preserves mass, energy and the divergence constraint of the magnetic field.
We then show some numerical results obtained with our approach.
We first test our scheme on simple academic test to compare the results with established methodologies, and then focus specifically on the simulation of plasma instabilities, with some tests on non Cartesian geometries to validate our discretization in the scope of tokamak instabilities.
\end{abstract}

\section{Introduction}
Over the last years, structure preserving discretizations have been receiving more and more attention \cite{gawlik2021variational,kraus2017gempic}. It is now understood that preserving some key properties such as Hamiltonian structure, symplecticity \cite{channell1990symplectic} or exterior derivative properties \cite{arnold2018finite} is fundamental in order to derive numerical schemes that behave well in long time. In several application such as geophysical fluid dynamics or plasma simulation in the context of inertial confinement fusion device, obtaining schemes that are able to produce accurate results in long time scale is of primal importance.

Several works have been conduced to derive structure preserving schemes from discrete least action principle, starting from finite dimensional systems, where the theory is well established (see \cite{marsden2001discrete} and reference therein), to more complicated continuum mechanics system, where the theory is still in development \cite{pavlov2011structure,gawlik2011geometric,natale2018variational,gawlik2020conservative,gawlik2021variational,gawlik2022finite}. 
For those infinite dimensional system, the approach is more complicated as one as to find good finite dimensional approximation for the system, while preserving the structure of the problem as much as possible.
One of the biggest limitation of the previously cited works on continuum systems lies in the fact that standard Hamilton principles are not able to encompass dissipative dynamics.
In order to overcome this issue, several techniques have been studied, such as the principle of least dissipation of energy, or minimum entropy production. Here we follow the recent approach of \cite{gay2017lagrangian1,gay2017lagrangian2} that uses a generalized Lagrange-d'Alembert principle. This approach was carried on to the discrete level in \cite{gawlik2022variational} were a discrete Lagrange-d'Alembert principle is stated for the viscous and heat conducting Euler equations. 

In this work we are interested in the viscoresistive equations of magnetohydrodynamics, which describe the evolution a magnetized plasma, taking in account two dissipative phenomena: viscosity, due to friction between fluid (plasma) particles, and resistivity, which models the imperfection of the plasma as a conductor. 
Those two non-ideal effects, are added to the variational principle of ideal magnetohydrodynamics following the ideas of \cite{gay2017lagrangian1,gay2017lagrangian2}. 
Then we use a similar procedure as the one presented in \cite{carlier2025variational} to derive our numerical scheme. 
The aim of this is to provide a stable discretization for the system of viscoresistive MHD, that is able to provide accurate results for long simulation specially for fusion devices.
Indeed traditional methods such as \cite{nikulsin2022jorek3d} suffer from the need of adding a lot of dissipation to stabilize the simulations and might therefore suppress some interesting physical behaviour.
The good results obtained for the ideal system in \cite{carlier2025variational} encourages us to follow the path of variational discretization to solve this issue.

The Lagrangian picture is known to be in duality with the Hamiltonian/Poisson bracket structure (one can go from one to another using a Legendre transform). 
When going to dissipative systems, several different theories have emerged to add the non-ideal effects to the symplectic structure. 
Among them \cite{kaufman1984dissipative} introduces the concept of dissipative bracket, 
\cite{morrison1984bracket} which provides an interpretation of the dissipation as a bracket with an entropy variable and
\cite{morrison1986paradigm} where this second bracket is interpreted as a metric term and introduces the so called "metriplectic systems".
Another formalism called "GENERIC", presented in \cite{grmela1997dynamics}, provides a dissipative bracket accounting for thermodynamic effect.
Recently \cite{morrison2024inclusive} proposes a unifying framework in which the dissipative effects are described by means of a 4-bracket, that has the same symmetries as a Riemann curvature tensor, leading to interesting geometric interpretations. 
In \cite{carlier2024metriplectic} we have shown that the generalized Lagrange-d'Alembert principle for dissipative system presented in \cite{gay2017lagrangian1,gay2017lagrangian2} can be in some cases rewrote using a metriplectic 4-bracket. In the present work we will be interested in this structure as we will show that this equivalence still holds at the semi-discrete level and will use the metriplectic form to derive a splitting time integration method. 

Motivated by the geometrical structure of the equations, we propose to use the Finite Element Exterior Calculus framework to discretize the different unknown present in the equation. 
This ensure the preservation of some key invariants such as total density and the solenoidal character of the magnetic field. Discretization is done using splines finite elements \cite{buffa2011isogeometric,back2012spline} as they are widely use for the discretization of fusion device. Those elements allow for for high order without the drawback of multiplying the number of degree of freedom and combine well with the use of coordinate aligned with equilibrium magnetic field. However the discretization described here does not rely on any particular property of spline elements and could be easily adapted to more standard FEEC spaces such as Nedelec or Raviart-Thomas \cite{raviart_mixed_1977,nedelec_mixed_1980,hiptmair2002finite}. 

The remainder of this article is organized as follow: in \cref{sec:vrmhd} we present the equation of viscoresistive MHD, the associated variational principle and the metriplectic reformulation. \Cref{sec:discrete} introduce the variational spatial discretization, its discrete metriplectic reformulation as well as the deduced time scheme, while \cref{sec:numerics} presents some numerical results obtained for the viscoresistive MHD equation as well as for the ideal MHD equation with some added artificial dissipation for stabilization. \Cref{sec:concl} will gather some concluding remarks.

\section{Variational Formulation for Viscous and Resistive Magnetohydrodynamics}
\label{sec:vrmhd}
\subsection{Viscous and resistive Magnetohydrodynamics}
The system of equations describing viscous and resistive magnetohydrodynamics (VRMHD) is given by \cite{goedbloed2004principles} : 
\begin{subequations}
\label{eqn:MHD}
\begin{equation}
\label{eqn:MHD_mom}
\rho \partial_t \uu + \rho (\uu \cdot \nabla \uu) + \nabla p + \BB \times \nabla \times \BB  = \nabla \cdot (\mu \nabla \uu) ~, 
\end{equation}
\begin{equation}
\label{eqn:MHD_mass}
\partial_t \rho + \nabla \cdot (\rho \uu) = 0 ~,
\end{equation}
\begin{equation}
\label{eqn:MHD_s}
\partial_t s + \nabla \cdot (s \uu) = \frac{1}{T}(\mu |\nabla \uu|^2 + \eta |\nabla \times \BB|^2 ) ~,
\end{equation}
\begin{equation}
\label{eqn:MHD_B}
\partial_t \BB + \nabla \times(\BB \times \uu) = - \nabla \times (\eta \nabla \times \BB) ~,
\end{equation}
\begin{equation}
\label{eqn:MHD_divB}
\nabla \cdot \BB = 0 ~,
\end{equation}
\end{subequations}
with $\rho$ the mass density, $\uu$ the velocity, $s$ the entropy and $B$ the magnetic field. 
$\mu$ is the viscosity of the considered plasma, $\eta$ its resistivity and $T$ denotes the temperature, defined here as a function of $s$ and $\rho$ (see next section).
\Cref{eqn:MHD_mom} describes the evolution of momentum of a particle of fluid under the different forces of pressure, Lorentz force and viscous friction with other particles.
The conservation of mass (continuity equation) is expressed by \cref{eqn:MHD_mass} and is standard. 
The equation for evolution of entropy \cref{eqn:MHD_s} is more involved, it describes the advection of entropy on the left hand side, while the right hand side is related to the increase of entropy caused by the non-conservative viscous and resistive forces.
\Cref{eqn:MHD_B} is the Faraday Law together with the Ohm law in a non-perfect conductor, it describe the evolution of the magnetic field in a resistive plasma.
We here make the choice of using the equation of the entropy (\cref{eqn:MHD_s}) while usual approaches use pressure of energy equations. This is in order to have the same variables in the variational principle below, which is an extension of the variational principle for ideal MHD, in which the entropy is purely advected (see for example \cite{carlier2025variational}.
\subsection{Variational formulation}
We now state a variational principle which solution are solution to the VRMHD equations. We will later use this principle to derive our numerical scheme, by stating a discrete variational principle and expressing its solution. The following formulation is inspired by the one for ideal MHD, using the framework described in \cite{gay2017lagrangian1,gay2017lagrangian2} to add the viscosity and a term that is new to our knowledge, also inspired by the development in \cite{gay2017lagrangian2} to include the resistivity.
\label{sec:var_form}
\begin{theorem}
\label{thm:var_form}
Consider the following Lagrangian :
\begin{equation}
l(\uu, \rho, s, \BB) = \int_{\Omega} \frac{\rho |\uu|^2}{2} - e(\rho, s) - \frac{|\BB|^2}{2}~,
\end{equation}
Where $e$ is the internal energy as a function of $\rho$ the density and $s$ the entropy. Solutions to the VRMHD momentum equation correspond to extremal curves of the action
\begin{equation}
S(\uu, \rho, s, \BB) = \int_0^T l(\uu(t), \rho(t), s(t), \BB(t)) dt
\end{equation}
Under the variational constraints :
\begin{equation}
\label{eqn:constraints}
\delta \uu = \partial_t \vv + [\uu, \vv], \qquad 
\delta \rho = - \nabla \cdot (\rho \vv)~, \qquad 
\frac{\delta l}{\delta s}(\delta s + \nabla \cdot(s \vv)) = -\mu \nabla \uu : \nabla \vv ~, 
\qquad \delta \BB = - \nabla \times (\BB \times \vv) ~.
\end{equation}
with $\vv$ a time dependent vector field that is null a time $t=0$ and $t=T$ and $[\uu,\vv] = \uu \cdot \nabla \vv - \vv \cdot \nabla \uu$.
\end{theorem}
To recover the viscous and resistive magnetohydrodynamic equations, this variational principle is supplemented with the following advection (with phenomenological constraint) equations :
\begin{subequations}
\begin{equation}
\label{eqn:advection_rho}
\partial_t \rho = - \nabla \cdot (\rho \uu) ~,
\end{equation}
\begin{equation}
\label{eqn:advection_s}
\frac{\delta l}{\delta s}(\partial_t s + \nabla \cdot(s \uu)) = -\mu |\nabla \uu|^2 - \eta |\nabla \times \BB|^2 ~,
\end{equation}
\begin{equation}
\label{eqn:advection_B}
\partial_t \BB = - \nabla \times (\BB \times \uu) - \nabla \times (\eta \nabla \times \BB) ~.
\end{equation}
\end{subequations}
The equation for the evolution of entropy is indeed equivalent to \cref{eqn:MHD_s} with the definition $T = - \frac{\partial l}{\partial s}$
\begin{proof}
Let $\uu$ be an extremal curve of the action $S = \int_0^T l(\uu, \rho, s, \BB)$ under the constraints \cref{eqn:constraints}.
Then for every curve $\vv$ in $X(\Omega)$ we have:
\begin{equation} \label{EL}
        0 = \frac{\delta S}{\delta \uu} \delta \uu + \frac{\delta S}{\delta \rho} \delta \rho + \frac{\delta S}{\delta s} \delta s + \frac{\delta S}{\delta \BB} \delta \BB 
        = \int_0^T \frac{\delta l}{\delta \uu} \delta \uu + \frac{\delta l}{\delta \rho} \delta \rho + \frac{\delta l}{\delta s} \delta s + \frac{\delta l}{\delta \BB} \delta \BB.
\end{equation}
%
%
Using integration by parts we develop the first term as
$$
\begin{aligned}
    \int_0^T \frac{\delta l}{\delta \uu} \delta \uu 
        &=  \int_0^T \int_\Omega \rho \uu \cdot (\partial_t \vv + \uu \cdot \nabla \vv - \vv \cdot \nabla \uu)
        \\
        &=    -\int_0^T \int_\Omega \partial_t(\rho \uu) \cdot \vv + \rho (\uu \cdot \nabla \uu) \cdot \vv + \uu \nabla \cdot(\rho \uu) \cdot \vv  + (\nabla \uu)^T (\rho \uu) \cdot \vv
\end{aligned}
$$
the second term as
$$
\begin{aligned}
    \int_0^T \frac{\delta l}{\delta \rho} \delta \rho
        &=  - \int_0^T \int_\Omega \Big( \frac{|\uu|^2}{2}-e(\rho, s)-\rho \partial_{\rho} e(\rho, s)\Big) \nabla \cdot(\rho \vv)
        \\
        &=  \int_0^T \int_\Omega \nabla \Big( \frac{|\uu|^2}{2} -e(\rho, s)-\rho \partial_{\rho} e(\rho, s)\Big) \cdot (\rho \vv)
        \\
        &=  \int_0^T \int_\Omega \big( (\nabla \uu)^T \uu - \nabla \rho \partial_{\rho} e(\rho, s) - \nabla s \partial_{s} e(\rho, s) - \nabla (\rho \partial_{\rho} e(\rho, s))\big) \cdot (\rho \vv)
\end{aligned}
$$
the third term as
$$
    \int_0^T \frac{\delta l}{\delta s} \delta s
        =  \int_0^T \int_\Omega \rho \partial_s e(\rho, s) \nabla \cdot (s \vv) - \mu \nabla \uu : \nabla \vv 
        =  \int_0^T \int_\Omega - \nabla (\rho \partial_s e(\rho, s)) \cdot (s \vv) + \nabla \cdot (\mu \nabla \uu) \cdot \vv 
$$
and the fourth term as
$$
    \int_0^T  \frac{\delta l}{\delta \BB} \delta \BB
        =  \int_0^T \int_\Omega  \BB \cdot \nabla \times(\BB \times \vv)
        =  -\int_0^T \int_\Omega  (\BB \times \nabla \times \BB) \cdot \vv ~.
$$
Introducing the pressure defined as $p=\rho(\rho \partial_{\rho} e+ s \partial_s e)$ we next observe that
$$
\rho \big(\nabla \rho \partial_{\rho} e(\rho, s) + \nabla s \partial_{s} e(\rho, s) + \nabla (\rho \partial_{\rho} e(\rho, s))\big) + s \nabla (\rho \partial_s e(\rho, s))
     = \nabla p~,
$$
so that \eqref{EL} being null for every $\vv$ yields

%
%
$$
0 
= \partial_t(\rho \uu) + \rho (\uu \cdot \nabla \uu) + \uu \nabla \cdot(\rho \uu) + \nabla p - \nabla \cdot (\mu \nabla \uu) + \BB \times \nabla \times \BB ~. 
$$
Finally, developing $\partial_t(\rho \uu) = \uu \partial_t \rho + \rho \partial_t \uu$ and using \cref{eqn:advection_rho} gives us the momentum equation for ideal MHD in its usual form: 
\begin{equation}
\label{eqn:momentum_MHD}
\rho \partial_t \uu + \rho (\uu \cdot \nabla \uu) + \nabla p - \nabla \cdot (\mu \nabla \uu) + \BB \times \nabla \times \BB  = 0~.
\end{equation}
\end{proof}
\subsection{Metriplectic reformulation}
\label{sec:metriplectic_continuous}
We now state another formulation that give the VRMHD equations, based on the metriplectic formalism. 
This framework is an extension of the symplectic framework in which the solution curves are described using two brackets (a standard poisson bracket for the non-dissipative part and a metric bracket for the dissipative one) with two different generators: the Hamiltonian corresponding to the total energy of the system and the entropy which is a Casimir (special invariant) of the non-dissipative system.
\begin{definition}[Metriplectic system]
\label{def:metriplectic}
Consider a symplectic 2-bracket $\{ \cdot , \cdot \}$ (that is antisymmetric, bilinear and satisfying the Jacobi identity),
and a metric 4-bracket $(\cdot,\cdot,\cdot,\cdot)$ with the following identity $(F,G;M,N) = -(G,F;M,N) = -(F,G;N,M) = (M,N;F,G)$ \cite{morrison2024inclusive}, both taking as argument functions of the dynamic.
Consider also two function of the dynamic : the Hamiltonian $H$ and the entropy $S$, that have the property that $\{F,S\}=0$ for all $F$.
The dynamic is a metriplectic system if for all function of the system $F$, we have $\dot{F} = \{F,H\} + (F,H;S,H)$
\end{definition}
\begin{remark}
If the dynamic satisfies $\dot{F} = \{F,H\}$ we have a standard Hamiltonian system, and the entropy is a Casimir of the system. The metric 4-bracket is responsible for the dissipation.
\end{remark}
In the case of MHD, we already know that symplectic bracket will be given by the standard Lie-Poisson bracket, the Hamiltonian should be the total energy and the entropy the total entropy. We only need to exhibit a metric 4-bracket that satisfies the identity above and generates the dissipative part when contracted with the Hamiltonian and the total entropy.
\begin{definition}
\label{def:mom_Ham}
We define the canonical momentum as 
\begin{equation}
\label{eqn:mom}
\mm = \frac{\partial l}{\partial \uu} = \rho \uu ~,
\end{equation}
the Hamiltonian (or energy) of the system as 
\begin{equation}
\label{eqn:ham}
H(\mm,\rho,s,\BB) = \langle \mm, \uu \rangle - l = \int_\Omega \frac{|\mm|^2}{2 \rho} + e(\rho,s) + \frac{|\BB|^2}{2}~,
\end{equation}
and the total entropy of the system 
\begin{equation}
\label{eqn:entr}
S(\mm,\rho,s,\BB) = \int_\Omega s ~.
\end{equation}
\end{definition}
We start by computing the variations of the Hamiltonian:
\begin{equation}
\frac{\delta H}{\delta \mm} = \uu ~, \qquad 
\frac{\delta H}{\delta \rho} = \frac{\partial \rho e}{\partial \rho} - \frac{|\uu|^2}{2} ~, \qquad 
\frac{\delta H}{\delta s} = \frac{\partial \rho e}{\partial s} ~, \qquad
\frac{\delta H}{\delta \BB} = \BB ~.
\end{equation}
Consider a curve solution to the VRMHD equation, we can rewrite the variational condition given by \cref{thm:var_form} as :
\begin{equation}
\int_\Omega \mm \cdot (\partial_t \vv + [\frac{\delta H}{\delta \mm}, \vv]) + \frac{\delta H}{\delta \rho} \nabla \cdot (\rho \vv) + \frac{\delta H}{\delta s} \nabla \cdot (s \vv) + \frac{\delta H}{\delta \BB} \nabla \times (\BB \times \vv) - \mu \nabla \frac{\delta H}{\delta \mm} \cdot \nabla \vv = 0~,
\end{equation}
Consider a function of the dynamic $F(\mm,\rho,s,\BB)$, we are interested in writing $\dot{F}$ with a symplectic bracket and a metric bracket. 
\begin{align*}
\dot{F} &= \int_\Omega \frac{\delta F}{\delta \mm} \dot{\mm} + \frac{\delta F}{\delta \rho} \dot{\rho} + \frac{\delta F}{\delta s} \dot{s} + \frac{\delta F}{\delta \BB} \dot{\BB} \\
		&= \int_\Omega \mm \cdot [\frac{\delta H}{\delta \mm},\frac{\delta H}{\delta \mm}] + \frac{\delta H}{\delta \rho} \nabla \cdot (\rho \frac{\delta F}{\delta \mm}) + \frac{\delta H}{\delta s} \nabla \cdot (s \frac{\delta F}{\delta \mm}) + \frac{\delta H}{\delta \BB} \cdot \nabla \times (\BB \times \frac{\delta F}{\delta \mm}) - \mu \nabla \frac{\delta H}{\delta \mm} \cdot \nabla \frac{\delta F}{\delta \mm} \\
		&- \frac{\delta F}{\delta \rho} \nabla \cdot (\rho \frac{\delta H}{\delta \mm}) - \frac{\delta F}{\delta s} \nabla \cdot (s \frac{\delta H}{\delta \mm}) + \frac{1}{T}\frac{\delta F}{\delta s}(\mu \nabla \frac{\delta H}{\delta \mm} \cdot \nabla \frac{\delta H}{\delta \mm} + \eta \nabla \times \frac{\delta H}{\delta \BB} \cdot \nabla \times \frac{\delta H}{\delta \BB}) \\
		&- \frac{\delta F}{\delta \BB} \cdot \nabla \times (\BB \times \frac{\delta H}{\delta \mm}) - \frac{\delta F}{\delta \BB} \nabla \times (\eta \nabla \times \frac{\delta H}{\delta \BB}) \\
		&= \{F,H\} + \int_\Omega \frac{1}{T}\frac{\delta F}{\delta s}(\mu \nabla \frac{\delta H}{\delta \mm} \cdot \nabla \frac{\delta H}{\delta \mm} + \eta \nabla \times \frac{\delta H}{\delta \BB} \cdot \nabla \times \frac{\delta H}{\delta \BB}) - \mu \nabla \frac{\delta H}{\delta \mm} \cdot \nabla \frac{\delta F}{\delta \mm} - \eta \nabla \times \frac{\delta F}{\delta \BB} \cdot \nabla \times \frac{\delta H}{\delta \BB}) \\
		&= \{F,H\} + \int_\Omega \frac{\mu}{T}(\frac{\delta F}{\delta s} \frac{\delta S}{\delta s} \nabla \frac{\delta H}{\delta \mm} \cdot \nabla \frac{\delta H}{\delta \mm} - \frac{\delta H}{\delta s} \frac{\delta S}{\delta s} \nabla \frac{\delta H}{\delta \mm} \cdot \nabla \frac{\delta F}{\delta \mm})\\
		& +\int_\Omega \frac{\eta}{T}(\frac{\delta F}{\delta s}\frac{\delta S}{\delta s}\nabla \times \frac{\delta H}{\delta \BB} \cdot \nabla \times \frac{\delta H}{\delta \BB} - \frac{\delta H}{\delta s}\frac{\delta S}{\delta s}\nabla \times \frac{\delta F}{\delta \BB} \cdot \nabla \times \frac{\delta H}{\delta \BB}) \\
		&= \{F,H\} + \int_\Omega \frac{\mu}{T}(\frac{\delta F}{\delta s} \frac{\delta S}{\delta s} \nabla \frac{\delta H}{\delta \mm} \cdot \nabla \frac{\delta H}{\delta \mm} - \frac{\delta H}{\delta s} \frac{\delta S}{\delta s} \nabla \frac{\delta H}{\delta \mm} \cdot \nabla \frac{\delta F}{\delta \mm} + \frac{\delta H}{\delta s} \frac{\delta H}{\delta s} \nabla \frac{\delta S}{\delta \mm} \cdot \nabla \frac{\delta F}{\delta \mm} - \frac{\delta F}{\delta s} \frac{\delta H}{\delta s} \nabla \frac{\delta S}{\delta \mm} \cdot \nabla \frac{\delta H}{\delta \mm})\\
		& +\int_\Omega \frac{\eta}{T}(\frac{\delta F}{\delta s}\frac{\delta S}{\delta s}\nabla \times \frac{\delta H}{\delta \BB} \cdot \nabla \times \frac{\delta H}{\delta \BB} - \frac{\delta H}{\delta s}\frac{\delta S}{\delta s}\nabla \times \frac{\delta F}{\delta \BB} \cdot \nabla \times \frac{\delta H}{\delta \BB} \\
		& + \frac{\delta H}{\delta s}\frac{\delta H}{\delta s}\nabla \times \frac{\delta S}{\delta \BB} \cdot \nabla \times \frac{\delta F}{\delta \BB} - \frac{\delta F}{\delta s}\frac{\delta H}{\delta s}\nabla \times \frac{\delta S}{\delta \BB} \cdot \nabla \times \frac{\delta H}{\delta \BB}) \\
		& = \{F,H\} + (F,H;S,H)_{visc} + (F,H;S,H)_{res}
\end{align*}
with 
\begin{subequations}
\label{eqn:brackets}
\begin{equation}
\begin{aligned}
\{F,G\} &= \int_\Omega \mm \cdot [\frac{\delta G}{\delta \mm},\frac{\delta H}{\delta \mm}] +  \frac{\delta G}{\delta \rho} \nabla \cdot (\rho \frac{\delta F}{\delta \mm}) + \frac{\delta G}{\delta s} \nabla \cdot (s \frac{\delta F}{\delta \mm}) + \frac{\delta G}{\delta \BB} \cdot \nabla \times (\BB \times \frac{\delta F}{\delta \mm}) \\
	&- \frac{\delta F}{\delta \rho} \nabla \cdot (\rho \frac{\delta G}{\delta \mm}) - \frac{\delta F}{\delta s} \nabla \cdot (s \frac{\delta G}{\delta \mm}) - \frac{\delta F}{\delta \BB} \cdot \nabla \times (\BB \times \frac{\delta G}{\delta \mm})
\end{aligned}
\end{equation}
\begin{equation}
(F,G,M,N)_{visc} = \int_\Omega \frac{\mu}{T}(\frac{\delta F}{\delta s} \frac{\delta M}{\delta s} \nabla \frac{\delta N}{\delta \mm} \cdot \nabla \frac{\delta G}{\delta \mm} - \frac{\delta G}{\delta s} \frac{\delta M}{\delta s} \nabla \frac{\delta N}{\delta \mm} \cdot \nabla \frac{\delta F}{\delta \mm} + \frac{\delta G}{\delta s} \frac{\delta N}{\delta s} \nabla \frac{\delta M}{\delta \mm} \cdot \nabla \frac{\delta F}{\delta \mm} - \frac{\delta F}{\delta s} \frac{\delta N}{\delta s} \nabla \frac{\delta M}{\delta \mm} \cdot \nabla \frac{\delta G}{\delta \mm})
\end{equation}
\begin{equation}
\begin{aligned}
(F,G,M,N)_{res} &=\int_\Omega \frac{\eta}{T} \Big(\frac{\delta F}{\delta s}\frac{\delta M}{\delta s}\nabla \times \frac{\delta N}{\delta \BB} \cdot \nabla \times \frac{\delta G}{\delta \BB} - \frac{\delta G}{\delta s}\frac{\delta M}{\delta s}\nabla \times \frac{\delta N}{\delta \BB} \cdot \nabla \times \frac{\delta F}{\delta \BB} \\
	& + \frac{\delta G}{\delta s}\frac{\delta N}{\delta s}\nabla \times \frac{\delta M}{\delta \BB} \cdot \nabla \times \frac{\delta F}{\delta \BB} - \frac{\delta F}{\delta s}\frac{\delta N}{\delta s}\nabla \times \frac{\delta M}{\delta \BB} \cdot \nabla \times \frac{\delta G}{\delta \BB} \Big)
\end{aligned}
\end{equation}
\end{subequations}
We have proved the following theorem
\begin{theorem}
A curve $(\uu(t),\rho(t), s(t), \BB(t))$ is solution to \cref{eqn:MHD} if and only if, doing the change of variable $(\uu, \rho, s, \BB) \mapsto (\mm, \rho, s, \BB)$, it is a metriplectic system \ref{def:metriplectic} with the brackets defined by \cref{eqn:brackets}, the Hamiltonian as in \cref{eqn:ham} and the entropy \cref{eqn:entr}.
\end{theorem}

\begin{remark}
Those 4-brackets have a Kulkarni-Nomizu product structure \cite{kulkarni1972bianchi,morrison2024inclusive}. This metriplectic system is similar to the one described in \cite{carlier2024metriplectic} and is found with similar computation. The major difference with the framework described in this work is the way resistivity is added, using additional terms on the advection equations for the entropy and magnetic field.
\end{remark}
\section{Discretization}
\label{sec:discrete}
\subsection{FEEC spaces}
Our discretization is based on the FEEC framework to discretize the different quantities. This allow to preserve crucial invariants, as it is based on discretizing the De Rham complex of differential forms.
The core property is to build spaces such that the calculus identities $\nabla \times \nabla = \nabla \cdot \nabla \times = 0$ hold. 

We consider four spaces $(V^0_h, V^1_h, V^2_h, V^3_h)$, such that:
\begin{itemize}
\item the gradient is well defined on $V_h^0$ ($V_h^0 \subset H^1(\Omega)$) and $\nabla(V_h^0) \subset V_h^1$.
\item the curl is well defined on $V_h^1$ ($V_h^1 \subset H(\text{curl},\Omega)$) and $\nabla \times (V_h^1) \subset V_h^2$.
\item the divergence is well defined on $V_h^2$ ($V_h^2 \subset H(\text{div},\Omega)$) and $\nabla \cdot(V_h^2) \subset V_h^3$.
\end{itemize}
Hence, we will look for a discrete density $\rhoh \in V_h^3$, entropy $\sh \in V_h^3$ and magnetic field $\Bh \in V_h^{2}$.

We also introduce the space $X_h=(V^0_h)^3$ that will be used to discretize the velocity $\uh \in X_h$

In the numerical simulations presented in \cref{sec:numerics}, we shall use spaces of tensor product splines of maximum regularity as described in \cite{buffa2011isogeometric}, however we point out that our method could be implemented with any other sequence of structure-preserving spaces satisfying this condition, such as Continuous Galerkin, Nédélec and Raviart-Thomas spaces \cite{hiptmair2002finite} or even sequences of Discontinuous Galerkin spaces \cite{gucclu2022broken}. Here we use the following spaces :

\begin{equation} \label{eqn:Discrete_spaces}
    \left\{ \begin{aligned}
    &V^0_h = S_{p+1} \otimes S_{p+1} \otimes S_{p+1} 
        \\
     &V^1_h = {  \left( \begin{smallmatrix}
     S_{p} &\otimes& S_{p+1} &\otimes& S_{p+1} \\ S_{p+1} &\otimes& S_{p} &\otimes& S_{p+1} \\ S_{p+1} &\otimes& S_{p+1} &\otimes& S_{p}
     \end{smallmatrix} \right) }      
     \\
     &V^2_h = 
     {  \left( \begin{smallmatrix}
     S_{p+1} &\otimes& S_{p} &\otimes& S_{p} \\ S_{p} &\otimes& S_{p+1} &\otimes& S_{p} \\ S_{p} &\otimes& S_{p} &\otimes& S_{p+1}
     \end{smallmatrix} \right)   } 
     \\
     &V^3_h = S_{p} \otimes S_{p} \otimes S_{p} ~      
    \end{aligned}  \right.
  \end{equation}
  
In order to write a discrete variational principle we will also need projections that maps from the continuous function space to the discrete ones. We therefore denote $\Pi^i$ a projection going from the continuous space to the discrete $V_h^i$. Here we use projectors of interpolation and histopolation based on geometric degrees of freedom
\cite{Bochev_Hyman_2006_csd,Robidoux.2008.histo,Gerritsma.2011.spec}.
\subsection{Discrete variational principle}
Our numerical scheme will be built using a discrete variational principle, mimicking the continuous one derived in \cref{sec:var_form}. We therefore consider the following discrete variational principle:
Consider the discrete Lagrangian
\begin{equation}
\label{eqn:discr_lag}
l_h(\uh,\rhoh,\sh,\Bh) = \int_\Omega \frac{1}{2}\rhoh |\uh|^2 - \rhoh e(\rhoh,sh) - \frac{1}{2}|\Bh|^2 ~,
\end{equation}
and the action 
\begin{equation}
\label{eqn:discr_action}
S_h = \int_0^T l_h dt
\end{equation}
\begin{problem}
\label{pbm:discr_var}
Find curves $\uh \in X_h, ~ \rhoh \in V^3_h, ~ \sh \in V^3_h, ~ \Bh \in V^2_h$ that extremize the action \cref{eqn:discr_action} under the constrained variations
\begin{subequations}
\label{eqn:discr_var}
\begin{equation}
\delta \uh = \partial_t \vh + \Pi^0([\uh,\vh]) ~, 
\end{equation}
\begin{equation}
\delta \rhoh = -\nabla \cdot \Pi^2(\rhoh \vh) ~, 
\end{equation}
\begin{equation}
\int_\Omega \frac{\delta l_h}{\delta \sh} (\delta \sh + \nabla \cdot \Pi^2(\rhoh \vh)) q_h = - \int \mu \nabla \uh : \nabla \vh q_h \forall q_h \in V^3_h ~,
\end{equation}
\begin{equation}
\delta \Bh = -\nabla \times \Pi^1(\Bh \times \vh) ~,
\end{equation}
\end{subequations}
with $\vh$ a curve in $X_h$ that is null at $t=0$ and $t=T$.

Together with this variational conditions, we have the following equations:
\begin{subequations}
\label{eqn:discr_adv}
\begin{equation}
\partial_t \rho_h = -\nabla \cdot \Pi^2(\rhoh \uh) ~,
\end{equation}
\begin{equation}
\partial_t \Bh = -\nabla \times \Pi^1(\Bh \times \uh) - \nabla \times (\eta \tilde{\nabla} \times \Bh) ~,
\end{equation}
\begin{equation}
\int_\Omega \frac{\delta l_h}{\delta \sh} (\partial_t \sh + \nabla \cdot \Pi^2(\rhoh \uh)) q_h = - \int (\mu |\nabla \uh|^2 + \eta |\tilde{\nabla} \times \Bh|^2) q_h \forall q_h \in V^3_h ~.
\end{equation}
\end{subequations}
\end{problem}
where we denote $\tilde{\nabla} \times$ the discrete dual of the curl that maps from $V^2_h$ to $V^1_h$ and defined by $\int_\Omega \tilde{\nabla} \times \Bh \cdot \Ah = \int_\Omega \Bh \cdot \nabla \times \Ah$ for $\Ah \in V^1_h$ and $\Bh \in V^2_h$

\subsection{Discrete FEM equations}
We now give the equations obtained by the previously stated discrete variational principle. The derivation is straightforward as we just plug the variational constraint and compute the variational derivative of the Lagrangian.
\begin{proposition}
The discrete equation corresponding to solution curves to \cref{pbm:discr_var} satisfy the following equations:
\begin{subequations}
\label{eqn:eq_semi_fem}
\begin{equation}
\label{eqn:semi_mom}
\begin{aligned}
\int_\Omega \partial_t (\rhoh \uh) \cdot \vh 
  - \rho_h \uh \cdot \Pi^0([\uh,\vh]) 
  + \Big(\frac{1}{2}|\uu_h|^2 - \frac{\partial \rho e}{\partial \rho}(\rho_h, s_h) \Big) \nabla \cdot\Pi^2(\rho_h \vh) \\
  - \frac{\partial \rho e}{\partial s}(\rho_h, s_h)  \nabla \cdot\Pi^2(s_h \vv_h) - \BB_h \cdot \nabla \times \Pi^1(\BB_h \times \vh) + \mu \nabla \uh : \nabla \vh &= 0 \qquad
  \forall \vv_h \in X_h
\end{aligned}
\end{equation}
\begin{equation}
\partial_t \rho_h + \nabla \cdot\Pi^2(\rho_h \uu_h)  = 0 ~, 
\end{equation}
\begin{equation}
\label{eqn:semi_entr}
\qquad \int_\Omega \frac{\partial e}{\partial s} (\partial_t \sh + \nabla \cdot \Pi^2(\rhoh \uh)) q_h = \int (\mu |\nabla \uh|^2 + \eta |\tilde{\nabla} \times \Bh|^2) q_h \qquad \forall q_h \in V^3_h ~, 
\end{equation}
\begin{equation}
\partial_t \Bh = -\nabla \times \Pi^1(\Bh \times \uh) - \nabla \times (\eta \tilde{\nabla} \times \Bh) ~.
\end{equation}
\end{subequations}
\end{proposition}
We introduce the weak weighting operators, defined for any function $f$ by 
\begin{equation}
M[f]_X : X_h \mapsto X_h ~, \qquad \int_\Omega M[f]_X \uh \cdot \vh = \int_\Omega f \uh \cdot \vh ~, \forall \uh, ~ \vh \in X_h ~,
\end{equation}
\begin{equation}
M[f]_3 : V^3_h \mapsto V^3_h ~, \qquad \int_\Omega M[f]_3 a_h b_h = \int_\Omega f a_h b_h ~, \forall a_h, ~ b_h \in V^3_h ~,
\end{equation}
and the projector $P_3$ which is the $L^2$ projector into $V^3_h$ defined by 
\begin{equation}
P_3 : F(\Omega) \mapsto V^3_h ~, \qquad \int_\Omega P_3(f) g_h = \int_\Omega f g_h ~, \forall g_h \in V^3_h ~.
\end{equation}
With those definition we can rewrite \cref{eqn:semi_mom} as 
\begin{equation}
\label{eqn:semi_mom_re}
\begin{aligned}
&\int_\Omega \partial_t (M[\rhoh]_X \uh) \cdot \vh 
  - M[\rhoh]_X \uh \cdot \big(\Pi^0([\uu_h, \vv_h])\big) 
  \\
  & \mspace{50mu}+ \Big(\frac{1}{2}|\uu_h|^2 - e(\rho_h, s_h) - \rho_h \partial_{\rho_h} e(\rho_h, s_h) \Big) \nabla \cdot\Pi^2(\rho_h \vh) \\
  & \mspace{50mu}
  - 	 \rho_h \partial_{s_h}e(\rho_h, s_h)  \nabla \cdot\Pi^2(s_h \vv_h) - \BB_h \cdot \nabla \times \Pi^1(\BB_h \times \vh) + \mu \nabla \uh : \nabla \vh = 0 \qquad
  \forall \vv_h \in X_h
\end{aligned}
\end{equation}
and \cref{eqn:semi_entr} 
\begin{equation}
\partial_t \sh + \nabla \cdot \Pi^2(\rhoh \uh) = M[T]_3^{-1} P_3(\mu |\nabla \uh|^2 + \eta |\tilde{\nabla} \times \Bh|^2)  ~.
\end{equation}

\subsection{Discrete metriplectic formulation}
We now show that the derived discrete equations also correspond to a discrete metriplectic system, that can be obtained in a very similar way as the continuous one was obtained in \cref{sec:metriplectic_continuous}. We start by defining discrete equivalent for the momentum, the Hamiltonian and the total entropy.
\begin{definition}
\label{def:mom_Ham_discr}
We define the discrete canonical momentum as 
\begin{equation}
\label{eqn:mom_discr}
\mh = \frac{\partial l_h}{\partial \uh} = M[\rhoh]_X \uh ~,
\end{equation}
the Hamiltonian (or energy) of the system as 
\begin{equation}
\label{eqn:ham_discr}
H_h(\mh,\rhoh,\sh,\Bh) = \langle \mh, \uh \rangle - l_h = \int_\Omega \frac{M[\rhoh]_X^{-1} \mh \cdot \mh}{2} + e(\rhoh,\sh) + \frac{|\Bh|^2}{2}~,
\end{equation}
and the total entropy of the system 
\begin{equation}
\label{eqn:entr_discr}
S_h(\mh,\rhoh,\sh,\Bh) = \int_\Omega \sh
\end{equation}
\end{definition}
We start by computing the variations of the Hamiltonian:
\begin{equation}
\frac{\delta H_h}{\delta \mh} = \uh ~, \qquad 
\frac{\delta H_h}{\delta \rhoh} = \frac{\partial \rho e}{\partial \rho} - \frac{|\uh|^2}{2} ~, \qquad 
\frac{\delta H_h}{\delta \sh} = \frac{\partial \rho e}{\partial s} ~, \qquad
\frac{\delta H_h}{\delta \Bh} = \Bh ~.
\end{equation}
We next remark that \cref{eqn:semi_mom_re} can be rewritten 
\begin{equation}
\begin{aligned}
&\int_\Omega \partial_t (\mh) \cdot \vh 
  - \mh \cdot \big(\Pi^0([\frac{\delta H_h}{\delta \mh}, \vv_h])\big) - \frac{\delta H_h}{\delta \rhoh} \nabla \cdot\Pi^2(\rho_h \vh) \\
  & \mspace{50mu}
  - 	 \frac{\delta H_h}{\delta \sh}  \nabla \cdot\Pi^2(s_h \vv_h) - \frac{\delta H_h}{\delta \Bh} \cdot \nabla \times \Pi^1(\BB_h \times \vh) + \mu \nabla \frac{\delta H_h}{\delta \mh} : \nabla \vh = 0 \qquad
  \forall \vv_h \in X_h ~.
\end{aligned}
\end{equation}
We now do a similar computation as in the continuous case, consider $F(\mh,\rhoh,\sh,\Bh)$ a function of the dynamic, we write:
\begin{align*}
\dot{F}&=\int \frac{\delta F}{\delta \mh} \dot{\mh} + \frac{\delta F}{\delta \rhoh} \dot{\rhoh} + \frac{\delta F}{\delta \sh} \dot{\sh} + \frac{\delta F}{\delta \Bh} \dot{\Bh} \\
	   &= \int_\Omega \mh \cdot \big(\Pi^0([\frac{\delta H_h}{\delta \mh}, \frac{\delta F}{\delta \mh}])\big) + \frac{\delta H_h}{\delta \rhoh} \nabla \cdot\Pi^2(\rho_h \frac{\delta F}{\delta \mh}) + \frac{\delta H_h}{\delta \sh}  \nabla \cdot\Pi^2(s_h \frac{\delta F}{\delta \mh}) \\
       & + \frac{\delta H_h}{\delta \Bh} \cdot \nabla \times \Pi^1(\BB_h \times \frac{\delta F}{\delta \mh}) - \mu \nabla \frac{\delta H_h}{\delta \mh} : \nabla \frac{\delta F}{\delta \mh} - \frac{\delta F}{\delta \rhoh}\nabla \cdot\Pi^2(\rhoh \frac{\delta H_h}{\delta \mh})\\
       & - \frac{\delta F}{\delta \sh} \nabla \cdot\Pi^2(\sh \frac{\delta H_h}{\delta \mh})  + \frac{\delta F}{\delta \sh} M[T]_3^{-1} P_3(\mu |\nabla \uh|^2 + \eta |\tilde{\nabla} \times \Bh|^2) - \frac{\delta F}{\delta \Bh} \cdot \nabla \times \Pi^1(\Bh \times \frac{\delta H_h}{\delta \mh})\\
       & - \frac{\delta F}{\delta \Bh}\nabla \times (\eta \tilde{\nabla} \times \frac{\delta H_h}{\delta \Bh}) \\
       & = \{F,H_h\}_h + \int_\Omega - \mu \nabla \frac{\delta H_h}{\delta \mh} : \nabla \frac{\delta F}{\delta \mh} + \frac{\delta F}{\delta \sh} M[T]_3^{-1} P_3(\mu |\nabla \uh|^2 + \eta |\tilde{\nabla} \times \Bh|^2) - \frac{\delta F}{\delta \Bh}\nabla \times (\eta \tilde{\nabla} \times \frac{\delta H_h}{\delta \Bh}) ~,
\end{align*}
with:
\begin{subequations}
\label{eqn:discr_simp_bracket}
\begin{equation}
\{F,G\}_h = \{F,G\}_{\mh} + \{F,G\}_{\rhoh} + \{F,G\}_{\sh} + \{F,G\}_{\Bh} ~,
\end{equation}
\begin{equation}
\{F,G\}_{\mh} = \int_\Omega \mh \cdot \big(\Pi^0([\frac{\delta G}{\delta \mh}, \frac{\delta F}{\delta \mh}])\big) ~,
\end{equation}
\begin{equation}
\{F,G\}_{\rhoh} = \int_\Omega \frac{\delta G}{\delta \rhoh} \nabla \cdot\Pi^2(\rhoh \frac{\delta F}{\delta \mh}) - \frac{\delta F}{\delta \rhoh} \nabla \cdot\Pi^2(\rhoh \frac{\delta G}{\delta \mh}) ~,
\end{equation}
\begin{equation}
\{F,G\}_{\sh} = \int_\Omega \frac{\delta G}{\delta \sh} \nabla \cdot\Pi^2(\sh \frac{\delta F}{\delta \mh}) - \frac{\delta F}{\delta \sh} \nabla \cdot\Pi^2(\sh \frac{\delta G}{\delta \mh}) ~,
\end{equation}
\begin{equation}
\{F,G\}_{\Bh} = \int_\Omega \frac{\delta G}{\delta \Bh} \cdot \nabla \times \Pi^1(\BB_h \times \frac{\delta F}{\delta \mh}) - \frac{\delta F}{\delta \Bh} \cdot \nabla \times \Pi^1(\BB_h \times \frac{\delta G}{\delta \mh}) ~.
\end{equation}
\end{subequations}
We now focus on the dissipative terms:
\begin{align*}
\dot{F} & = \{F,H_h\}_h + \int_\Omega - \mu \nabla \frac{\delta H_h}{\delta \mh} : \nabla \frac{\delta F}{\delta \mh} + \frac{\delta F}{\delta \sh} M[T]_3^{-1} P_3(\mu |\nabla \uh|^2 + \eta |\tilde{\nabla} \times \Bh|^2) - \frac{\delta F}{\delta \Bh}\nabla \times (\eta \tilde{\nabla} \times \frac{\delta H_h}{\delta \Bh}) \\
	    & = \{F,H_h\}_h + \int_\Omega \frac{\delta S_h}{\delta \sh} \frac{\delta F}{\delta \sh} M[T]_3^{-1} P_3(\mu \nabla \frac{\delta H_h}{\delta \mh}:\nabla \frac{\delta H_h}{\delta \mh}) - \frac{\delta S_h}{\delta \sh} M[T]_3 M[T]_3^{-1} P_3(\mu \nabla \frac{\delta H_h}{\delta \mh} : \nabla \frac{\delta F}{\delta \mh}) \\
	    & + \int_\Omega \frac{\delta S_h}{\delta \sh} \frac{\delta F}{\delta \sh} M[T]_3^{-1} P_3(\eta \tilde{\nabla} \times \frac{\delta H_h}{\delta \Bh} \cdot \tilde{\nabla} \times \frac{\delta H_h}{\delta \Bh}) - \frac{\delta S_h}{\delta \sh} M[T]_3 M[T]_3^{-1} P_3(\eta \tilde{\nabla} \times \frac{\delta H_h}{\delta \Bh} \cdot \tilde{\nabla} \times \frac{\delta F}{\delta \Bh}) \\
	    & = \{F,H_h\} + (F,H_h;S_h,H_h)_{h,visc} + (F,H_h;S_h,H_h)_{h,res} ~,
\end{align*}
with
\begin{subequations}
\label{eqn:discr_metric_bracket}
\begin{equation}
\begin{aligned}
(F,G,M,N)_{h,visc} &= \int_\Omega (\frac{\delta F}{\delta \sh} \frac{\delta M}{\delta \sh} M[T]^{-1}P_3(\mu \nabla \frac{\delta N}{\delta \mh} \cdot \nabla \frac{\delta G}{\delta \mh}) - \frac{\delta G}{\delta \sh} \frac{\delta M}{\delta \sh} M[T]^{-1}P_3(\mu \nabla \frac{\delta N}{\delta \mh} \cdot \nabla \frac{\delta F}{\delta \mh})\\
   & + \frac{\delta G}{\delta \sh} \frac{\delta N}{\delta \sh} M[T]^{-1}P_3(\mu \nabla \frac{\delta M}{\delta \mh} \cdot \nabla \frac{\delta F}{\delta \mh}) - \frac{\delta F}{\delta \sh} \frac{\delta N}{\delta \sh} M[T]^{-1}P_3(\mu \nabla \frac{\delta M}{\delta \mh} \cdot \nabla \frac{\delta G}{\delta \mh}) ~,
\end{aligned}
\end{equation}
\begin{equation}
\begin{aligned}
(F,G,M,N)_{h,res} &= \int_\Omega \frac{\delta F}{\delta \sh}\frac{\delta M}{\delta \sh} M[T]^{-1}P_3(\eta \nabla \times \frac{\delta N}{\delta \Bh} \cdot \nabla \times \frac{\delta G}{\delta \Bh}) - \frac{\delta G}{\delta \sh}\frac{\delta M}{\delta \sh} M[T]^{-1}P_3(\eta \nabla \times \frac{\delta N}{\delta \Bh} \cdot \nabla \times \frac{\delta F}{\delta \Bh}) \\
   & + \frac{\delta G}{\delta \sh}\frac{\delta N}{\delta \sh} M[T]^{-1}P_3(\eta \nabla \times \frac{\delta M}{\delta \Bh} \cdot \nabla \times \frac{\delta F}{\delta \Bh}) - \frac{\delta F}{\delta \sh}\frac{\delta N}{\delta \sh} M[T]^{-1}P_3(\eta \nabla \times \frac{\delta M}{\delta \Bh} \cdot \nabla \times \frac{\delta G}{\delta \Bh}) ~.
\end{aligned}
\end{equation}
\end{subequations}
All the brackets obtained are consistent discretization of the continuous ones derived in \cref{sec:metriplectic_continuous}, however we can see that they do not correspond to the first naive discretization of the bracket one could do.
\begin{theorem}
A curve $(\uh(t),\rhoh(t), \sh(t), \Bh(t))$ is solution to \cref{pbm:discr_var} if and only if, doing the change of variable $(\uh, \rhoh, \sh, \Bh) \mapsto (\mh, \rhoh, \sh, \Bh)$, it is a metriplectic system \ref{def:metriplectic} with the brackets defined by \cref{eqn:discr_simp_bracket,eqn:discr_metric_bracket}, the Hamiltonian as in \cref{eqn:ham_discr} and the entropy \cref{eqn:entr_discr}.
\end{theorem}
\subsection{Bracket splitting and time discretization}
Using the previous results, we clearly see in \cref{eqn:discr_simp_bracket} that the symplectic bracket can be split in several smaller brackets that all have the entropy as a casimir. 
This fact encourage us to use bracket splitting method as time discretization. 
This type of time discretization integrates the systems using several subsystems that are each defined by a small bracket. 
Solving a smaller subsystem allows to define an \textit{integrator} that can then be composed in various ways to obtain high order time integration. 
We here explicit the different subsystems and respective integrator that are given by every small bracket, by $a^n$ we denote a quantity at the $n$-th time step, $\Delta t$ denotes the time step and $a^{n+\frac{1}{2}} = \frac{a^{n}+a^{n+1}}{2}$.

\underline{$\{F,G\}_{\mh}$ :}

Evolve only $\uh$ and is given by:
\begin{equation}
\int_\Omega \rhoh^n \frac{\uh^{n+1}-\uh^n}{\Delta t}\cdot \vh - \int_\Omega \rhoh^n \uh^n \cdot \Pi^0([\uh^{n+\frac{1}{2}}, \vh]) = 0 \qquad \forall \vh \in X_h ~.
\end{equation} 
We denote the corresponding integrator by $\Phi^{\mh}_{\Delta t}$.

\underline{$\{F,G\}_{\rhoh}$ :}

Evolve $\uh$ and $\rhoh$:
\begin{subequations}
\begin{equation}
\int_\Omega  \frac{\rhoh^{n+1} \uh^{n+1}-\rhoh^{n} \uh^n}{\Delta t}\cdot \vh + \int_\Omega \Big(\frac{\uh^{n+1}\cdot\uh^n}{2} - \frac{\rhoh^{n+1}e(\rhoh^{n+1},\sh^n)-\rhoh^{n}e(\rhoh^{n},\sh^n)}{\rhoh^{n+1}-\rhoh^n}\Big)\nabla \cdot (\Pi^2(\rhoh^n \vh)) = 0 \qquad \forall \vh \in X_h ~,
\end{equation} 
\begin{equation}
\label{eqn:adv_rhoh}
\frac{\rhoh^{n+1}-\rhoh^n}{\Delta t} + \nabla \cdot \Pi^2(\rhoh^n \uh^{n+\frac{1}{2}}) = 0 ~.
\end{equation}
\end{subequations}
We denote the corresponding integrator by $\Phi^{\rhoh}_{\Delta t}$.

\underline{$\{F,G\}_{\sh}$ :}

Evolve $\uh$ and $\sh$:
\begin{subequations}
\begin{equation}
\int_\Omega \rhoh^{n} \frac{\uh^{n+1}- \uh^n}{\Delta t}\cdot \vh - \int_\Omega \Big(\frac{\rhoh^{n}e(\rhoh^n,\sh^{n+1})-\rhoh^{n}e(\rhoh^{n},\sh^n)}{\sh^{n+1}-\sh^n}\Big)\nabla \cdot (\Pi^2(\sh^n \vh)) = 0 \qquad \forall \vh \in X_h ~,
\end{equation} 
\begin{equation}
\label{eqn:adv_sh}
\frac{\sh^{n+1}-\sh^n}{\Delta t} + \nabla \cdot \Pi^2(\sh^n \uh^{n+\frac{1}{2}}) = 0 ~.
\end{equation}
\end{subequations}
We denote the corresponding integrator by $\Phi^{\sh}_{\Delta t}$.

\underline{$\{F,G\}_{\Bh}$ :}

Evolve $\uh$ and $\Bh$:
\begin{subequations}
\begin{equation}
\int_\Omega \rhoh^{n} \frac{\uh^{n+1}- \uh^n}{\Delta t}\cdot \vh - \int_\Omega \Bh^n \cdot \nabla \times (\Pi^1(\Bh^{n} \times \vh)) = 0 \qquad \forall \vh \in X_h ~,
\end{equation} 
\begin{equation}
\label{eqn:adv_Bh}
\frac{\Bh^{n+1}-\Bh^n}{\Delta t} + \nabla \times \Pi^1(\Bh^n \times \uh^{n+\frac{1}{2}}) = 0 ~.
\end{equation}
\end{subequations}
We denote the corresponding integrator by $\Phi^{\Bh}_{\Delta t}$.

\underline{$(F,G,M,N)_{h,visc}$ :}

Evolve $\uh$ and $\sh$:
\begin{subequations}
\begin{equation}
\int_\Omega \rhoh^n \frac{\uh^{n+1}- \uh^n}{\Delta t}\cdot \vh + \int_\Omega \mu \nabla \uh^{n+1} : \nabla \vh = 0 \qquad \forall \vh \in X_h ~,
\end{equation} 
\begin{equation}
\int_\Omega \frac{\rhoh^n e(\rhoh^n,\sh^{n+1})-\rhoh^n e(\rhoh^n,\sh^n)}{\Delta t} q_h - \int_\Omega \mu \nabla \uh^{n+\frac{1}{2}} : \nabla \uh^{n+1} q_h = 0 \qquad \forall q_h \in V^3_h ~.
\end{equation}
\end{subequations}
We denote the corresponding integrator by $\Phi^{visc}_{\Delta t}$.

\underline{$(F,G,M,N)_{h,res}$ :}

Evolve $\Bh$ and $\sh$:
\begin{subequations}
\begin{equation}
\frac{\Bh^{n+1}- \Bh^n}{\Delta t} + \nabla \times(\eta \tilde{\nabla} \times \Bh^{n+1}) = 0  ~,
\end{equation} 
\begin{equation}
\int_\Omega \frac{\rhoh^n e(\rhoh^n,\sh^{n+1})-\rhoh^n e(\rhoh^n,\sh^n)}{\Delta t} q_h - \int_\Omega \eta \nabla \times \Bh^{n+\frac{1}{2}} \cdot \nabla \times \Bh^{n+1} q_h = 0 \qquad \forall q_h \in V^3_h ~.
\end{equation}
\end{subequations}
We denote the corresponding integrator by $\Phi^{res}_{\Delta t}$.
\begin{proposition}
All the integrator previously defined preserve the total mass $\int_\Omega \rhoh$, solenoidal character of $\BB$ ($\nabla \cdot \Bh$ = 0) and total energy $H_h$. The integrator corresponding to symplectic parts ($\Phi^{\mh}_{\Delta t}$, $\Phi^{\rhoh}_{\Delta t}$, $\Phi^{\sh}_{\Delta t}$ and $\Phi^{\Bh}_{\Delta t}$) all preserve the total entropy $S_h$. 
\end{proposition}
\begin{remark}
For the two dissipative propagators ($\Phi^{visc}_{\Delta t}$ and $\Phi^{res}_{\Delta t}$) we opted for an implicit time discretization of the dissipation term, in order to avoid restrictive CFL-like condition on the time step. This prevent us from formally proving that the create entropy, however it is the case as long as the variable at time $n+1$ is "not too far" from the one at time $n$.
\end{remark}
In numerical experiments, we use the Strang splitting (guaranteeing higher order time integration), that is the integrator for a full time step is given by
\begin{equation}
\label{eqn:Strang_splitting}
\Phi_{\Delta t} = \Phi^{\rhoh}_{\Delta t/2} \circ \Phi^{\mh}_{\Delta t/2} \circ \Phi^{\sh}_{\Delta t/2} \circ \Phi^{\Bh}_{\Delta t/2} \circ \Phi^{visc}_{\Delta t/2} \circ \Phi^{res}_{\Delta t/2} \circ \Phi^{res}_{\Delta t/2} \circ \Phi^{visc}_{\Delta t/2} \circ \Phi^{\Bh}_{\Delta t/2} \circ \Phi^{\sh}_{\Delta t/2} \circ \Phi^{\mh}_{\Delta t/2} \circ \Phi^{\rhoh}_{\Delta t/2} ~.
\end{equation}
\begin{proposition}
The scheme defined by the total integrator \cref{eqn:Strang_splitting} preserves the total mass, energy and solenoidal character of $\BB$.
\end{proposition}
\subsection{Artificial viscosity and resistivity}
One of the possibility offered by our framework is the use of artificial viscosity to stabilize simulations of ideal MHD. 
For more references about artificial viscosity we refer to \cite{gentry1966eulerian,harlow1971numerical}, here we use first derivative of the interest quantities to scale the dissipative term.

In the previous section, the viscosity parameter $\mu$ and the resistive parameter $\eta$ were taken as constant, however the whole framework described here is also valid for non constant parameters.
We here propose to use a variable parameter in order to only add dissipation in strong gradient zones. In this case we use 
\begin{equation}
\mu = \mu_a |\nabla \uh|_2 ~,
\end{equation}
\begin{equation}
\eta = \eta_a |\nabla \times \Bh|_2 ~,
\end{equation}
where $\eta_a$ and $\mu_a$ are user defined parameters, usually scaled with the mesh size. In our case we use $\mu_a = \eta_a = 2 h^2$.
This type of artificial dissipation originate in \cite{vonneumann1950method} and is very basic and more involved methodology could be used (for example using shock detector to only add dissipation where needed). However our goal here is only to present and test this possibility and not to conduct a study on the different types of artificial dissipation.
\section{Numerical Examples}
\label{sec:numerics}
We now present some numerical results obtained with this scheme. The main goal here is to present the improvements from the scheme presented in \cite{carlier2025variational}, that is in one hand the ability to stabilize ideal simulation, presented \cref{sec:ideal_test}, and in the other hand the integration of viscous and resistive test in \cref{sec:vr_test}

In all our numerical experiments the internal energy of the plasma is given by $e(\rho, s) = \rho^{\gamma-1} \exp(s/\rho)$ which is a standard equation of state for perfect gaz. 
\subsection{Stabilized ideal tests}
\label{sec:ideal_test}
\subsubsection{Dispersion relation study}
\begin{figure}
\centering
\captionsetup[subfigure]{justification=centering}
    \begin{subfigure}[b]{0.3\textwidth}
    \centering
    \includegraphics[width=\textwidth]{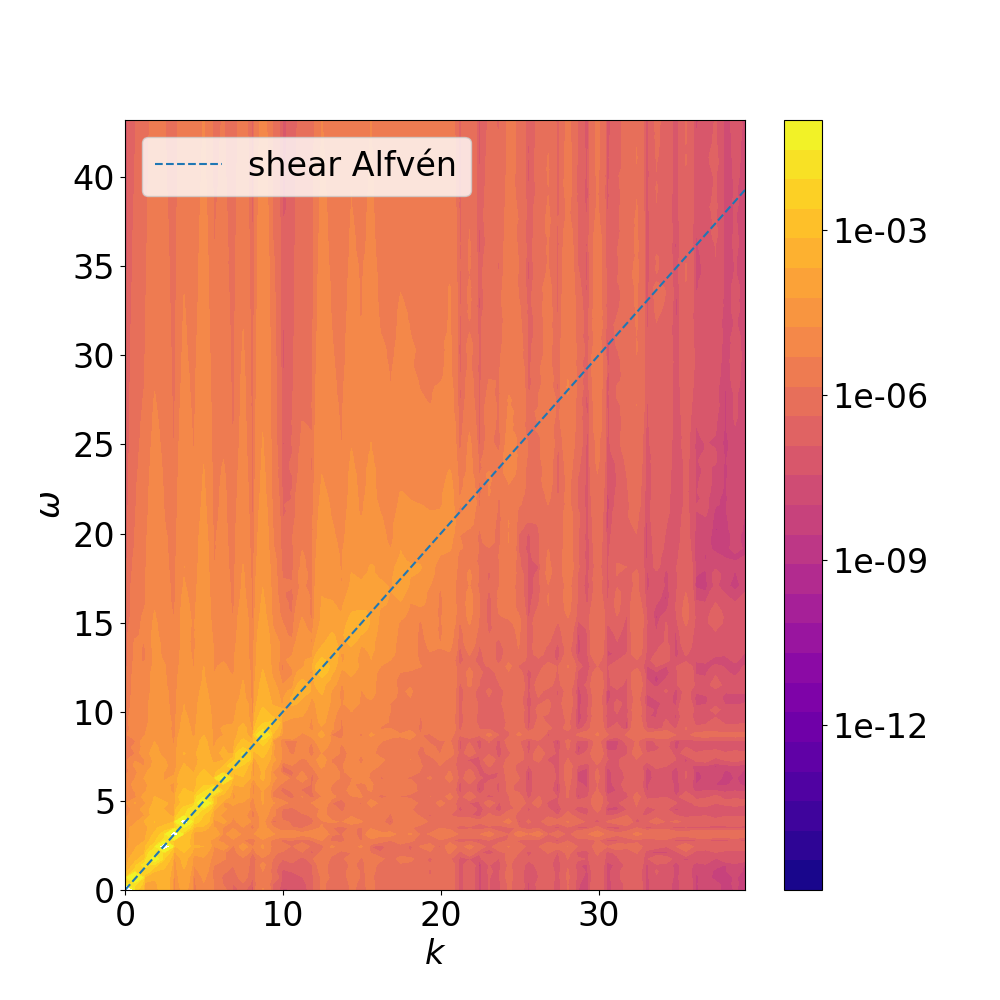}
    \caption{$\uu$ power spectrum for a $128$ elements grid with artificial dissipation}
    \end{subfigure}
    \begin{subfigure}[b]{0.3\textwidth}
    \centering
    \includegraphics[width=\textwidth]{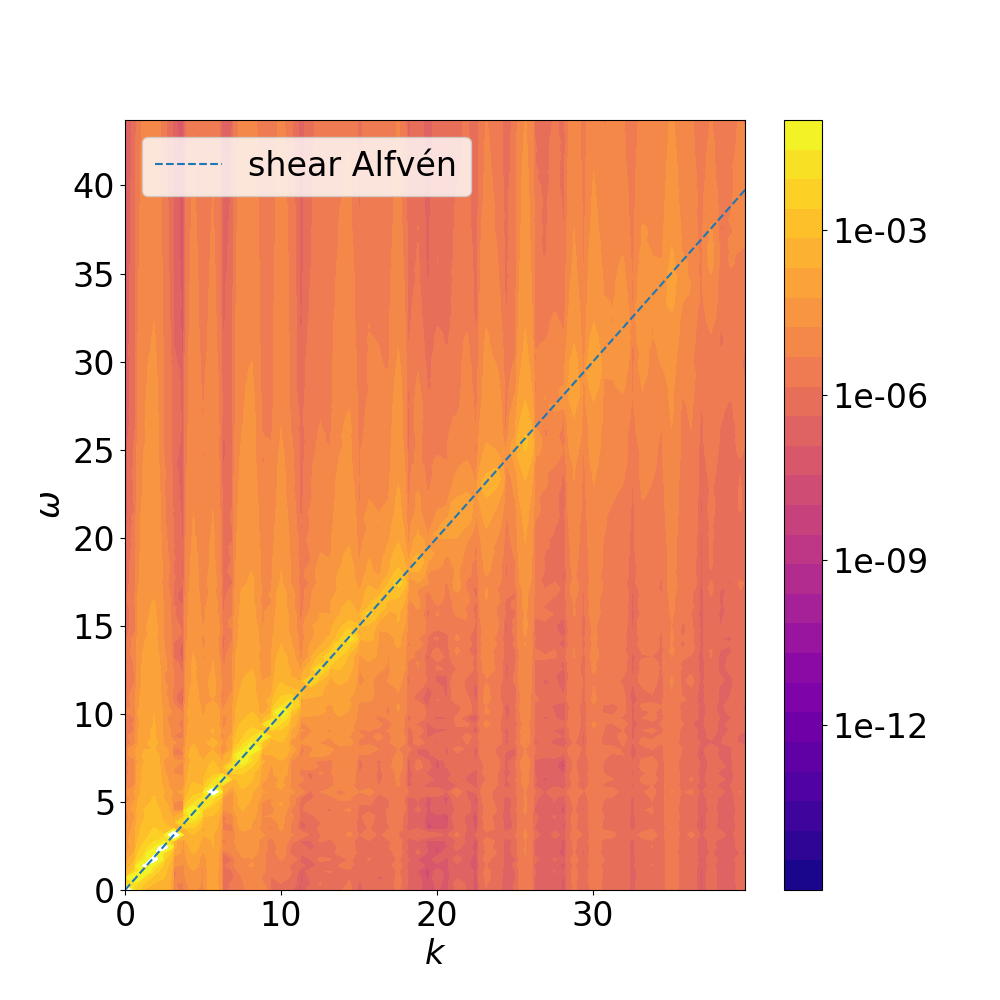}
    \caption{$\uu$ power spectrum for a $256$ elements grid with artificial dissipation}
    \end{subfigure}
    \begin{subfigure}[b]{0.3\textwidth}
    \centering
    \includegraphics[width=\textwidth]{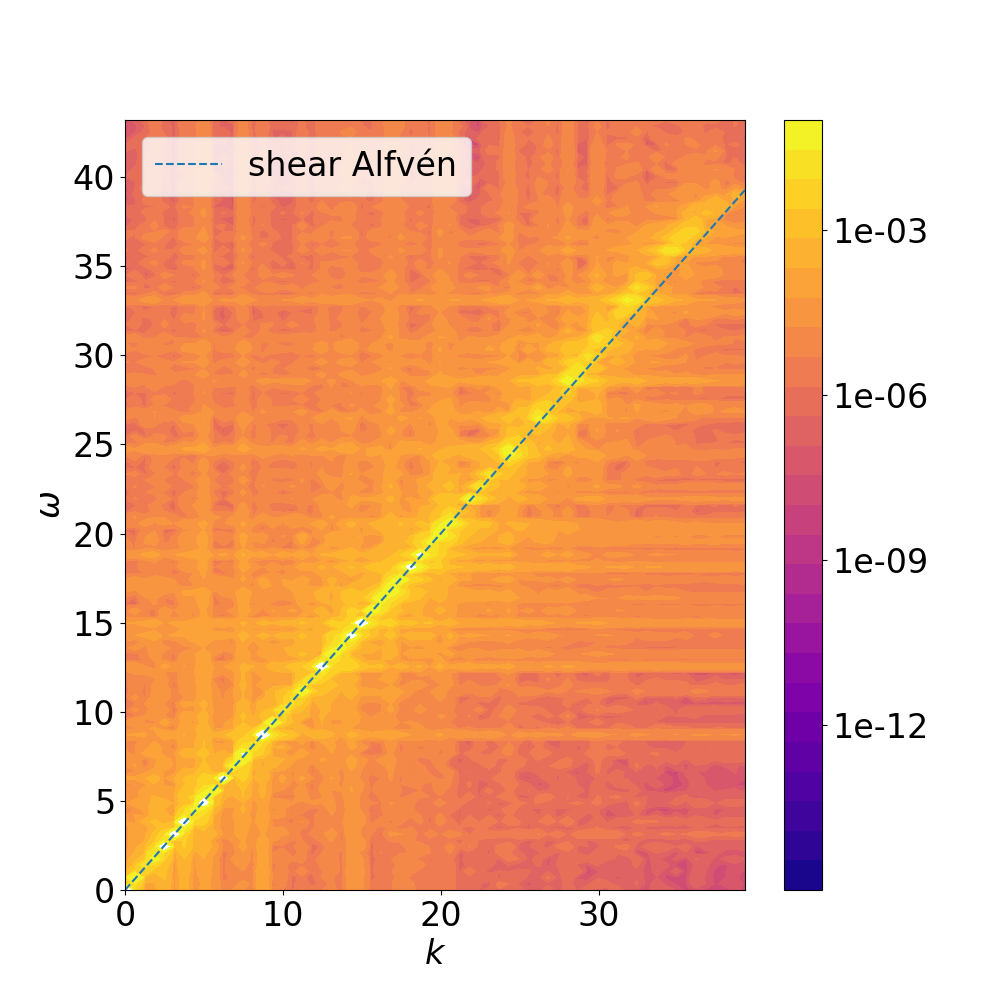}
    \caption{$\uu$ power spectrum for a $128$ elements grid without artificial dissipation}
    \end{subfigure}
    \begin{subfigure}[b]{0.3\textwidth}
    \centering
    \includegraphics[width=\textwidth]{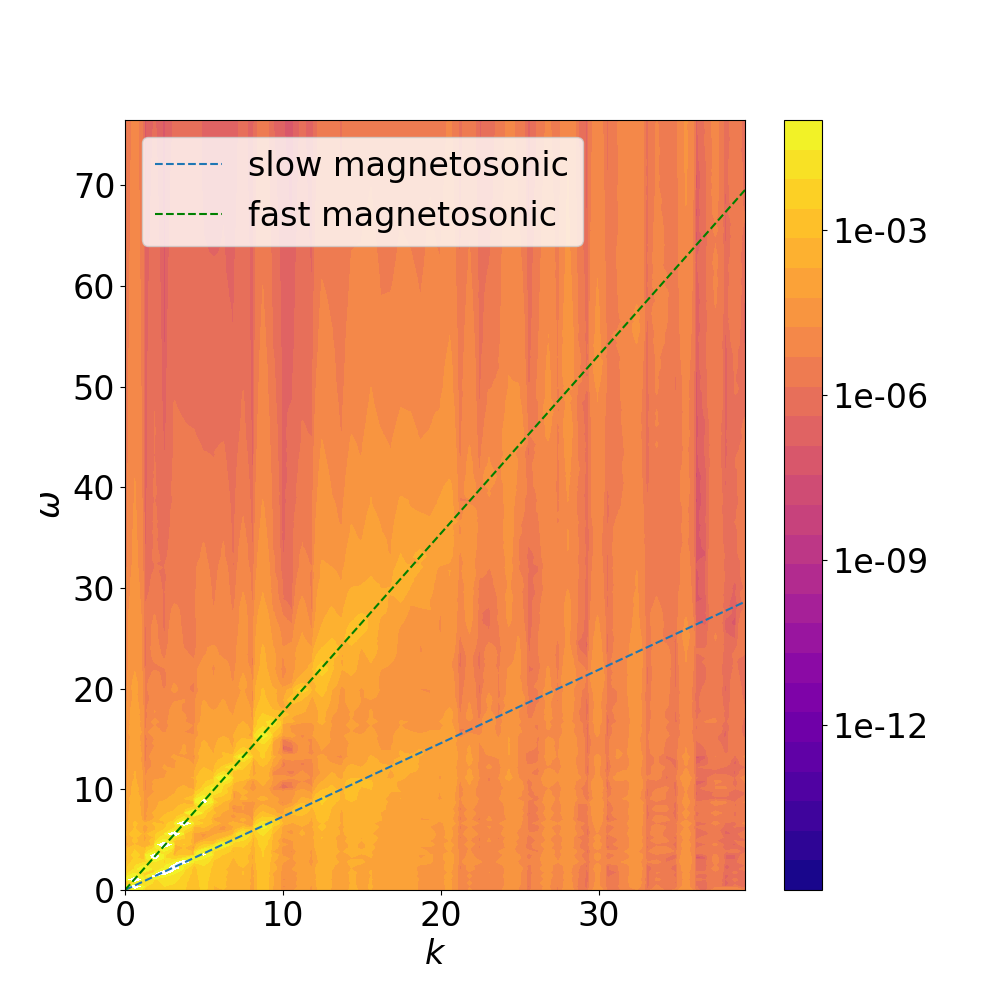}
    \caption{$p$ power spectrum for a $128$ elements grid with artificial dissipation}
    \end{subfigure}
    \begin{subfigure}[b]{0.3\textwidth}
    \centering
    \includegraphics[width=\textwidth]{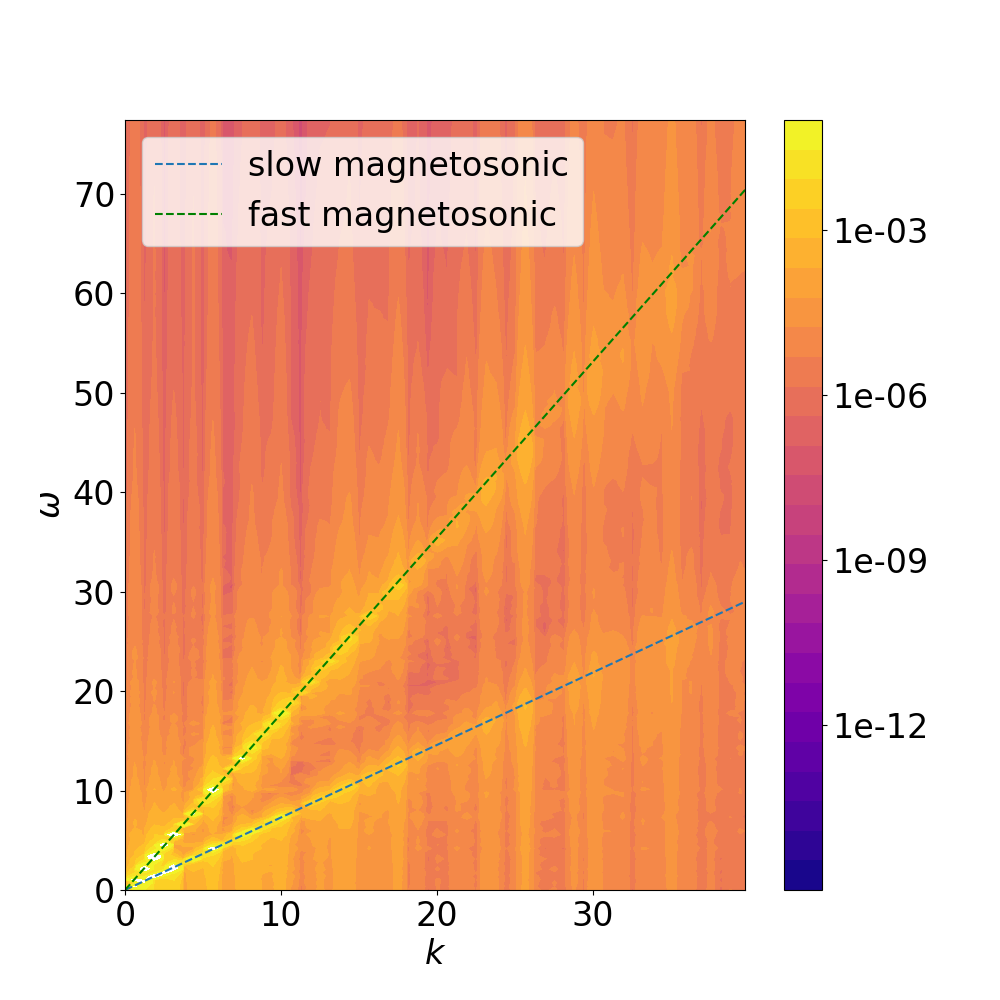}
    \caption{$p$ power spectrum for a $256$ elements grid with artificial dissipation}
    \end{subfigure}
    \begin{subfigure}[b]{0.3\textwidth}
    \centering
    \includegraphics[width=\textwidth]{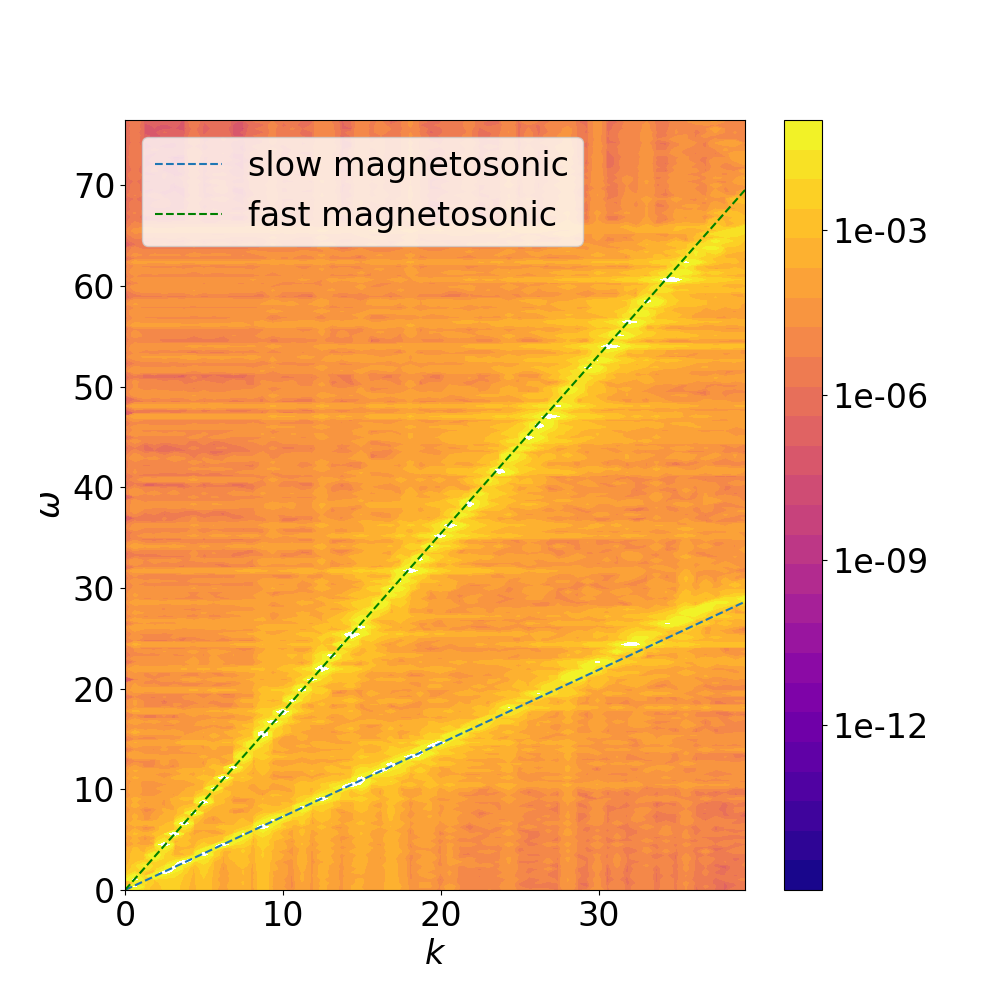}
    \caption{$p$ power spectrum for a $128$ elements grid without artificial dissipation}
    \end{subfigure}
\caption{Evolution of the dispersion relation with the use of artificial dissipation}
\label{fig:disp_rel}
\end{figure}

Our first test is a dispersion relation study, trying to reproduce and compared with the results obtained on \cite{holderied2021mhd}. We study the propagation of waves in the $x-$ direction (that is waves that can be written as $\exp(i(kx-\omega t)$). For such wave in an homogeneous magnetic field $\BB = B_x \mathbf{e}_x + B_y \mathbf{e}_y$, the dispersion relation reads
\begin{equation}
\label{eqn:dispersion_relation}
\Big(\omega^2 - k^2 v_A^2 \frac{B_x^2}{B_x^2 + B_y^2} \Big) \Big[\omega^2-\frac{1}{2}k^2(c_s^2 + v_A^2)(1 \pm \sqrt{\delta}) \Big] = 0 ~, \qquad \delta = \frac{4 B_x^2 c_s^2 v_A^2}{(c_s^2 + v_A^2)^2(B_x^2 + B_y^2)} ~,
\end{equation}
with 
\begin{equation}
v_A^2 = \frac{B_x^2 + B_y^2}{\rho} ~ \qquad c_S^2 = \gamma \frac{p}{\rho}
\end{equation}
To study the evolution of this wave with our scheme, we set a one dimensional domain $[0,10]$, initialized with constant density $\rho=1$, magnetic field $\BB = (1,1,0)$ and entropy $s = \log(\frac{1}{(\gamma-1)})$. We then set the initial velocity as random noise with amplitude $10^{-2}$, to excite the whole spectrum of frequency.  
In \cref{eqn:dispersion_relation} the term in the parenthesis correspond to shear Alven wave that will be found in the spectrum of $\uu$, while the bracket term correspond to slow and fast magnetosonic waves in the spectrum of the pressure. 
The results are presented in \cref{fig:disp_rel}. They are obtained after running the simulation until $t=18.$ with a time step of $\Delta t = 3 \times 10^{-2}$. The two leftmost ones show stabilized simulation with $128$ and $256$ elements in the $x$ direction, using stabilization parameters $\mu_a = \eta_a = 2 h^2 \approx 1.2 \times 10^{-2}$ for the coarser one and $3.0 \times 10^{-3}$ for the finer one. The right one shows the results for an unstabilized simulation, with $128$ element. All simulation where run with $p=2$ as maximum spline degree. In all cases we see that the dispersion is well solved and not subject to deviation as the linear model presented in \cite{holderied2021mhd}, thus showing a great improvement. However we see that with numerical dissipation the higher frequency modes are damped and do not carry as much energy as the lower ones. This behaviour is normal, as dissipation tends to damp faster the high modes. We also observe that with a finer grid, as the artificial dissipation goes smaller, the higher modes are better resolved.

\subsubsection{Ideal Orszag-Tang Vortex}
Our second test is an ideal Orszag-Tang vortex, the domain is a two dimensional periodic square $[0, 2 \pi]^2$ and the simulation is setup with the following initial conditions:
\begin{align*}
\rho(x,y,0) &= \gamma^2 ~, \\
s(x,y,0) &= \gamma^2 \log\Big(\frac{\gamma}{(\gamma-1)\gamma^{2\gamma}}\Big) ~, \\
\uu (x,y,0) &= (-\sin(y), \sin(x)) ~,\\
\BB (x,y,0) &= (-\sin(y), \sin(2x)) ~ 
\end{align*}
with $\gamma = 5/3$. We use a grid of $256 \times 256$ elements, with a spline degree of 2. A constant time step $\Delta t = 10^{-3}$ is used and artificial dissipation parameters are set to $\mu_a = \eta_a = 2h^2 \approx 1.2 \times 10^{-3}$. This tests exhibit a shock soon before $t=1.$ .
\begin{figure}
\centering
\captionsetup[subfigure]{justification=centering}
    \begin{subfigure}[b]{0.49\textwidth}
    \centering
    \includegraphics[width=\textwidth]{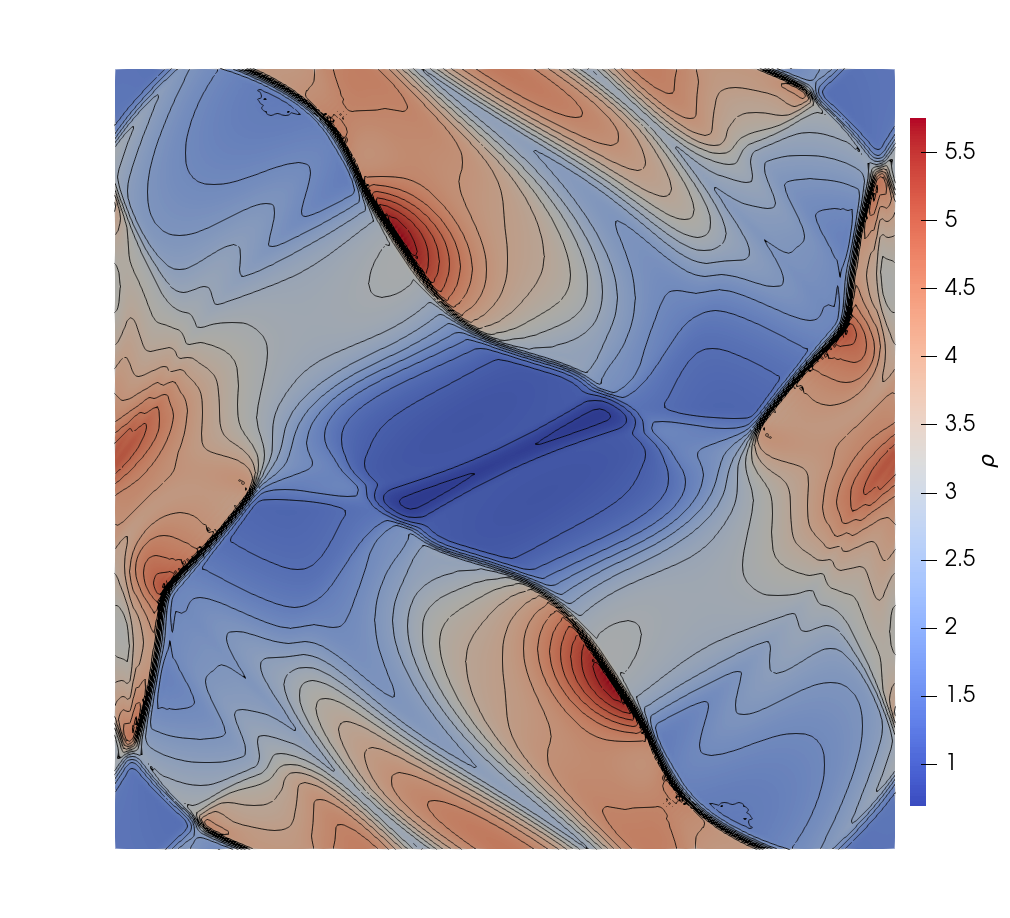}
    \caption{Density with equispaced contour every $0.25$}
    \end{subfigure}
    \begin{subfigure}[b]{0.49\textwidth}
    \centering
    \includegraphics[width=\textwidth]{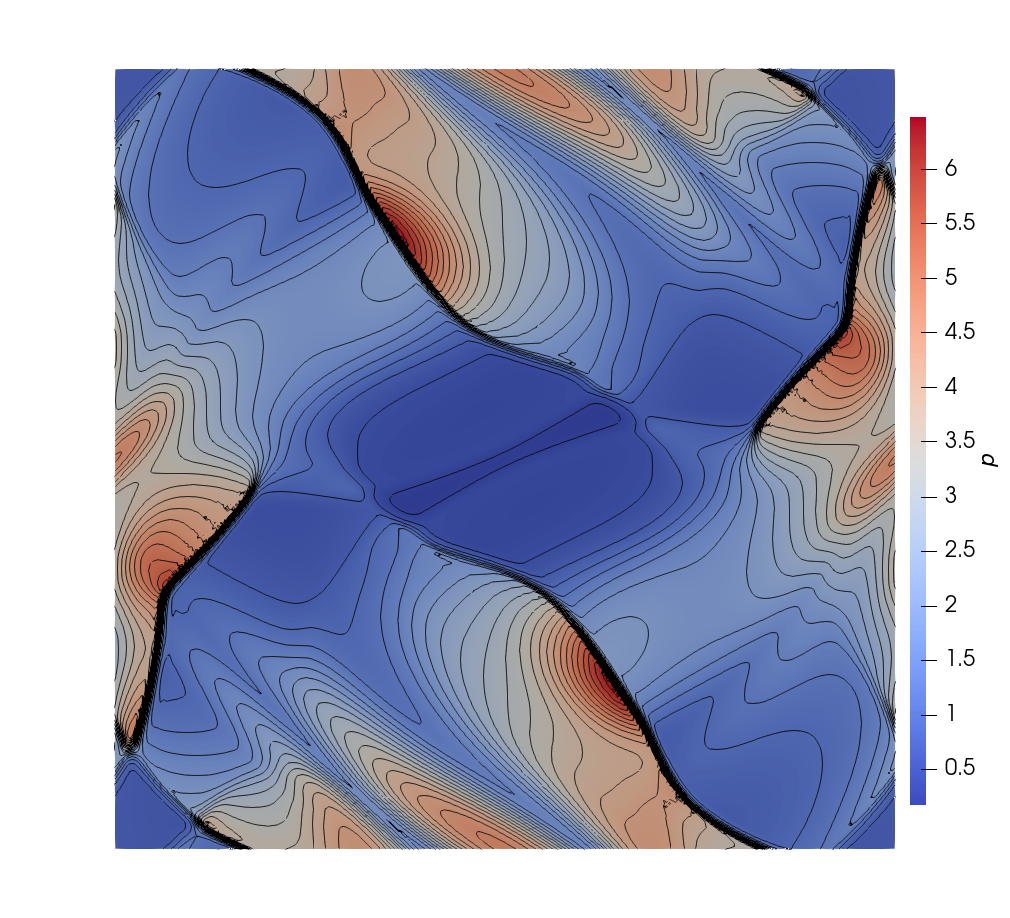}
    \caption{Pressure with equispaced contour every $0.25$}
    \end{subfigure}
    \begin{subfigure}[b]{0.49\textwidth}
    \centering
    \includegraphics[width=\textwidth]{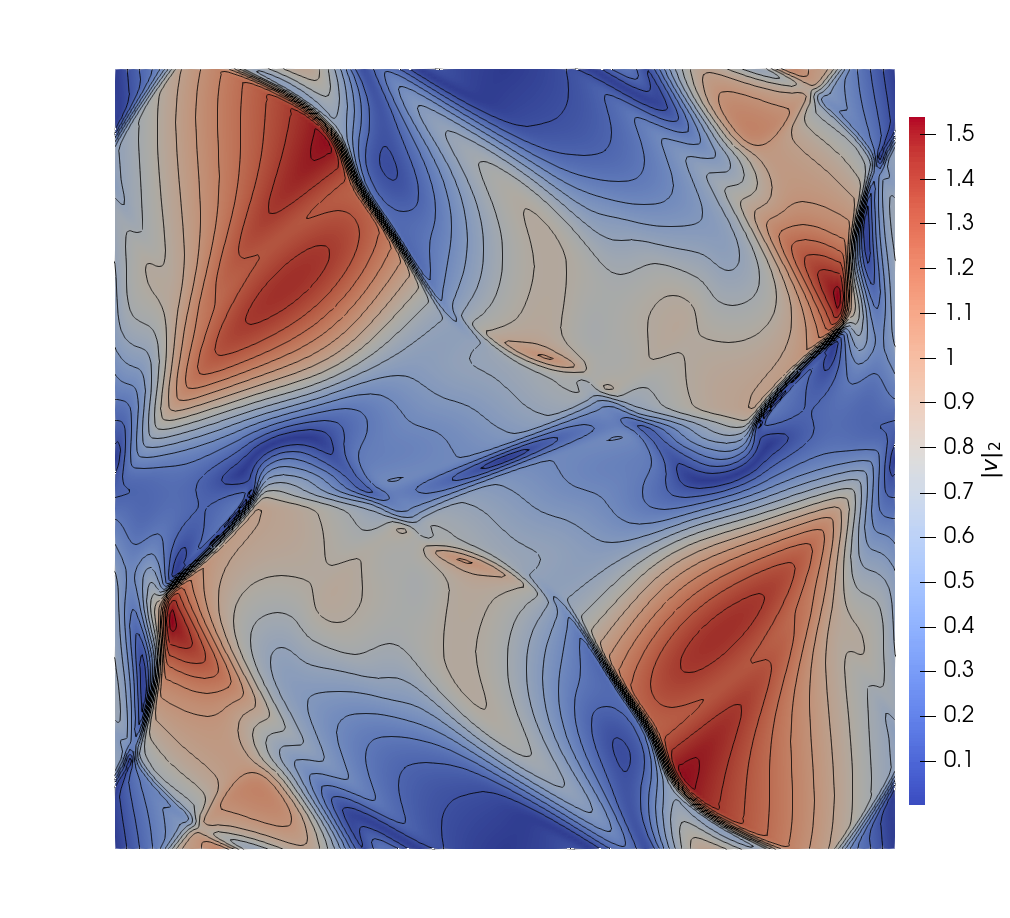}
    \caption{$L^2$ norm of the velocity with equispaced contour every $0.1$}
    \end{subfigure}
    \begin{subfigure}[b]{0.49\textwidth}
    \centering
    \includegraphics[width=\textwidth]{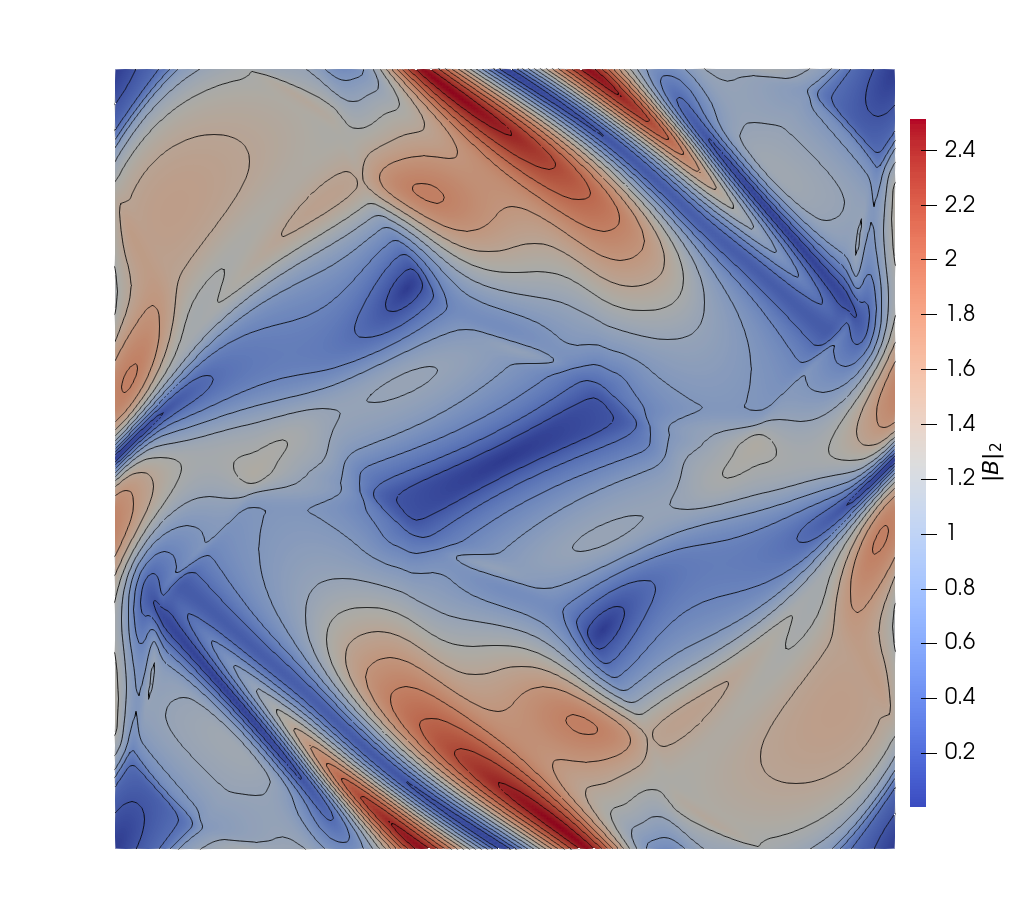}
    \caption{$L^2$ norm of the magnetic field with equispaced contour every $0.2$}
    \end{subfigure}
\caption{Orszag-Tang vortex at $t=2$}
\label{fig:OT}
\end{figure}
\begin{figure}
\centering
\captionsetup[subfigure]{justification=centering}
    \begin{subfigure}[b]{0.24\textwidth}
    \centering
    \includegraphics[width=\textwidth]{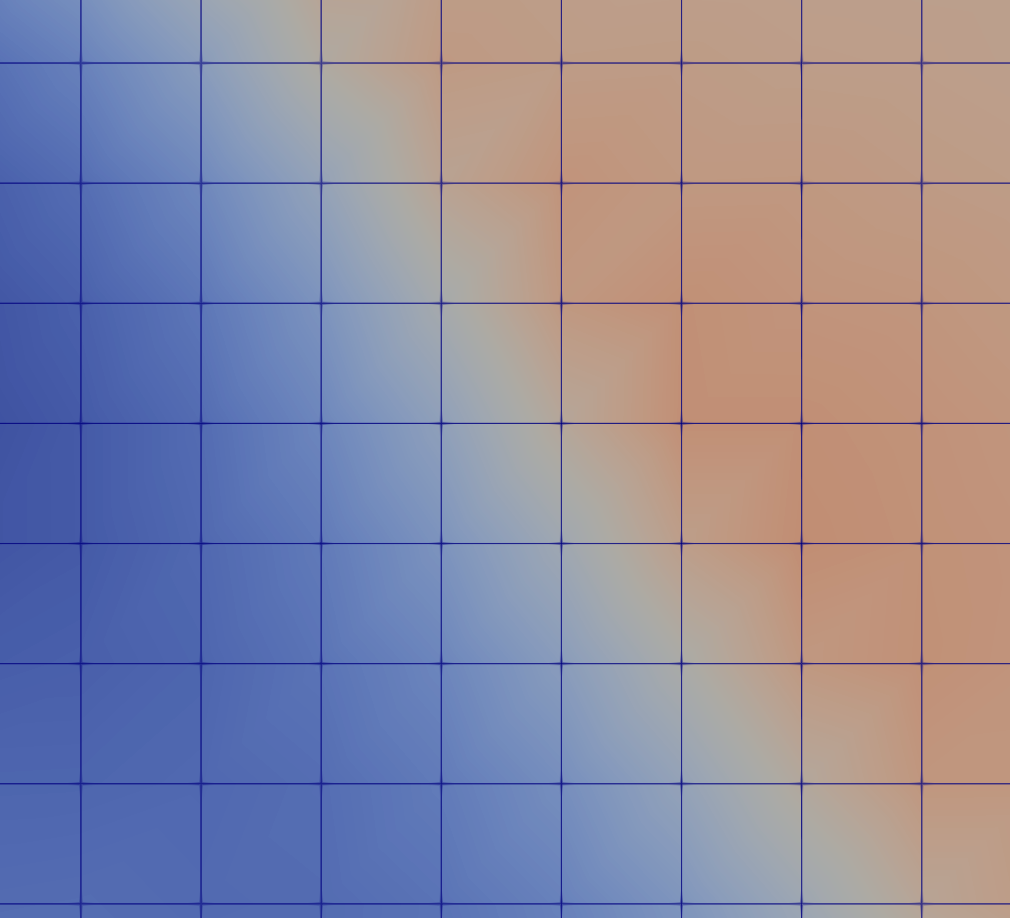}
    \caption{$64\times 64$ grid}
    \end{subfigure}
    \begin{subfigure}[b]{0.24\textwidth}
    \centering
    \includegraphics[width=\textwidth]{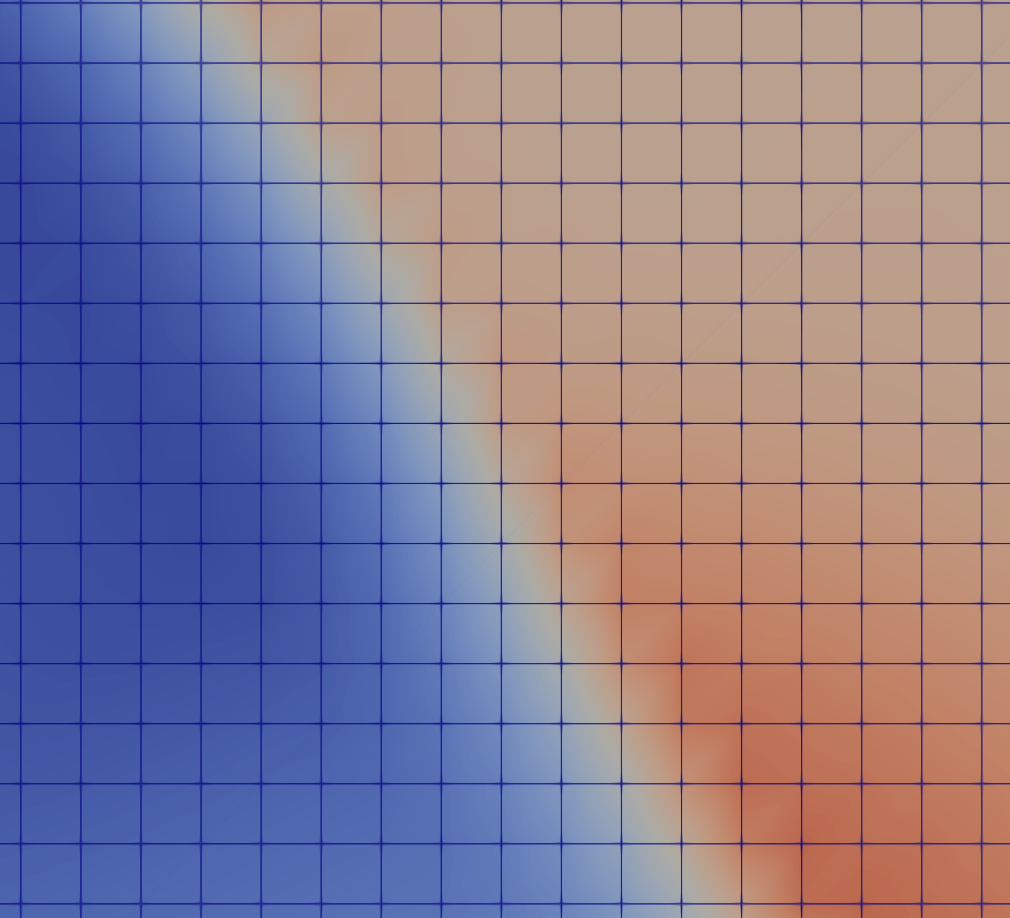}
    \caption{$128\times 128$ grid}
    \end{subfigure}
    \begin{subfigure}[b]{0.24\textwidth}
    \centering
    \includegraphics[width=\textwidth]{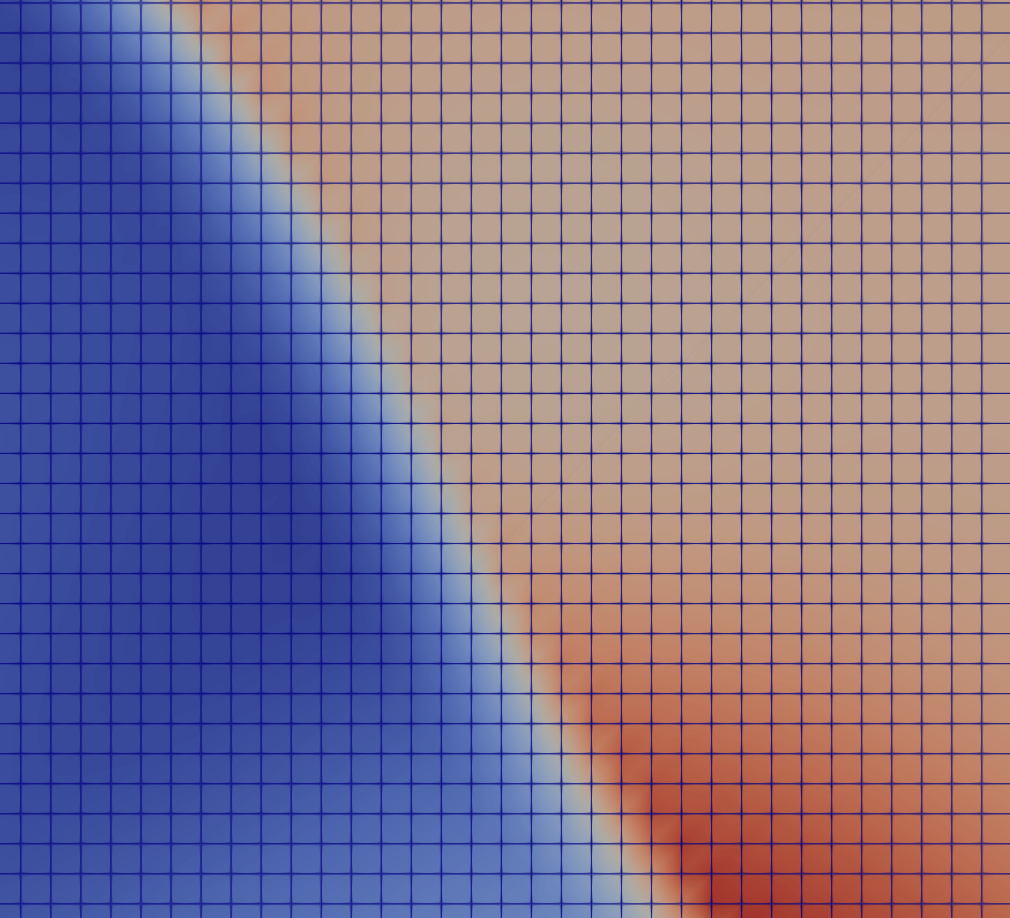}
    \caption{$256\times 256$ grid}
    \end{subfigure}
    \begin{subfigure}[b]{0.24\textwidth}
    \centering
    \includegraphics[width=\textwidth]{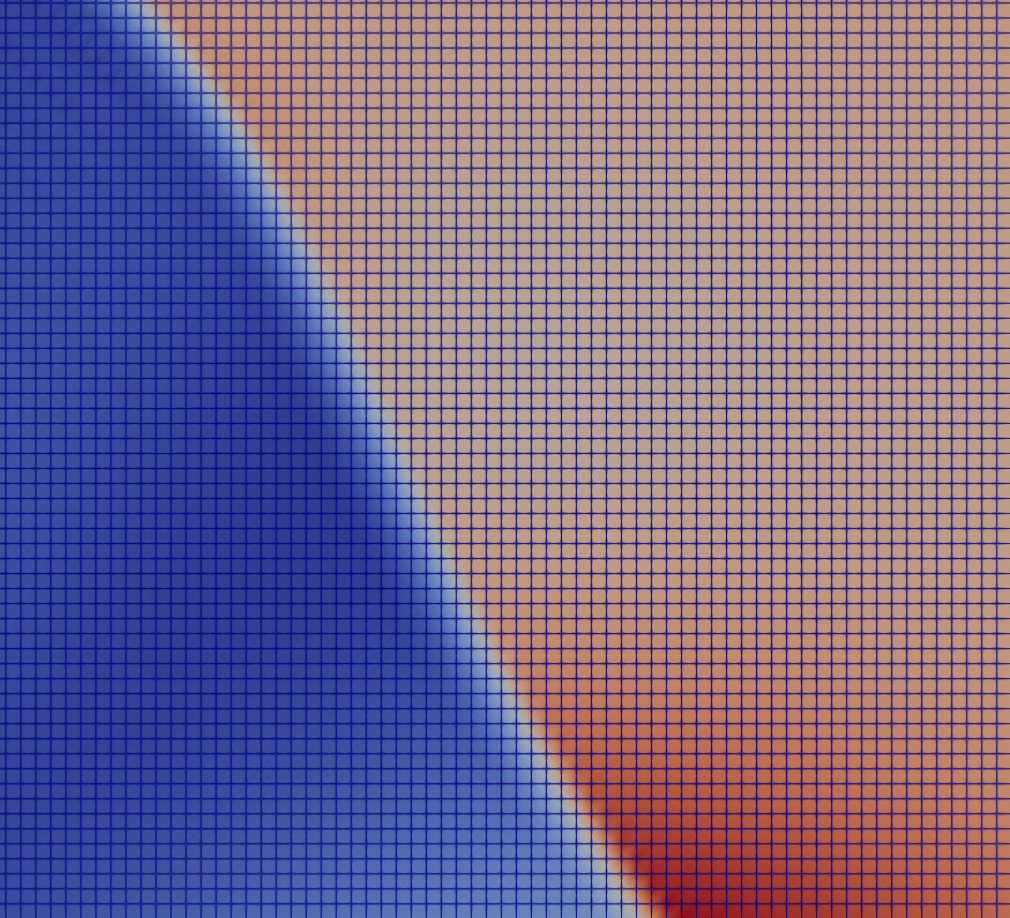}
    \caption{$512\times 512$ grid}
    \end{subfigure}
\caption{Orszag-Tang vortex at $t=2$, zoom on the shock with different meshes}
\label{fig:OT_shock}
\end{figure}
The results are presented in \cref{fig:OT}. We can see that our variational scheme is now able to reproduce well the dynamic of this test. With the addition of artificial viscosity and resistivity we don't observe any spurious oscillation and see that the shocks are correctly captured. 
As expected they suffer from some dissipation, but are resolved within 3 mesh elements, independently on the mesh size, as shown in \cref{fig:OT_shock} showing good ability of our scheme to handle discontinuities.
We also point that our scheme is able to compute accurately the pressure, although it's not a primary variable (it has to be calculated from the density and entropy), without suffering from oscillations.
\subsubsection{Ideal Kelvin-Helmholtz instability}
We present another application of the stabilization using artificial viscosity, a Kelvin-Helmholtz instability.
The setup for this test is a periodic rectangle domain $[0,1] \times [0,2]$ with initial conditions
\begin{align*}
\rho(x,y,0) &= 0.5 + 0.75 T_\delta(y) ~, \\
s(x,y,0) &= -\rho(x,y,0)(\log(\gamma-1)+\gamma \log(\rho(x,y,0))) ~, \\
u_x (x,y,0) &= 0.5(T_\delta(y)-1) ~,\\
u_y (x,y,0) &= 0.1 \sin(2 \pi x) ~, \\
\text{with } T_\delta(y) &= -\tanh((y - 0.5)/\delta)+\tanh((y + 0.5)/\delta) ~,
\end{align*}
where $\delta=\frac{1}{15}$. The entropy is set up so that initially the pressure is constant over the domain. This test is purely fluid ($\BB = 0$).
\begin{figure}
\centering
\captionsetup[subfigure]{justification=centering}
    \begin{subfigure}[b]{0.49\textwidth}
    \centering
    \includegraphics[width=\textwidth]{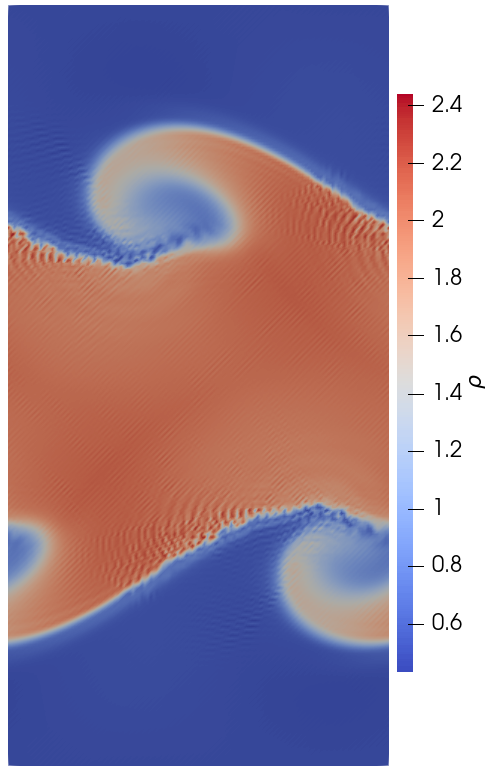}
    \caption{Without artificial viscosity}
    \end{subfigure}
    \begin{subfigure}[b]{0.49\textwidth}
    \centering
    \includegraphics[width=\textwidth]{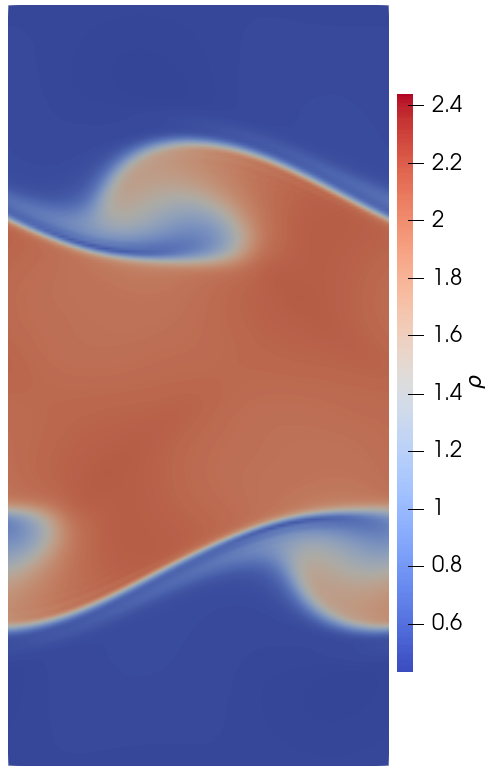}
    \caption{With artificial viscosity}
    \end{subfigure}
\caption{Density for the Kelvin-Helmholtz instability at $t=2$}
\label{fig:khi}
\end{figure}
The results are presented in \cref{fig:khi}. On the left panel we present the result of the simulation without using artificial viscosity while the right panel shows results using stabilization. Both simulations are run on a $128 \times 256$ grid, with maximal degree of spline $p=2$, constant time step of $\Delta t = 5 \times 10^{-4}$ and $\gamma = \frac{7}{5}$. For the dissipative simulation we keep the choice $\mu_a = 2 \times h^2 \approx 1.2 10^{-4}$. Results show that artificial viscosity is indeed able to stabilize the simulation, as on the non-stabilized one, we observe oscillations (that would lead the code to crash due to negative density if ran for too long), while the stabilized simulation doesn't show any sign of instability. We point out that the triggered Kelvin-Helmholtz instability is damped by the dissipation.
\subsubsection{Toroidal Alfven Eigenwave}
\label{sec:tae_test}
\begin{figure}
    \centering
    \includegraphics[width=\textwidth]{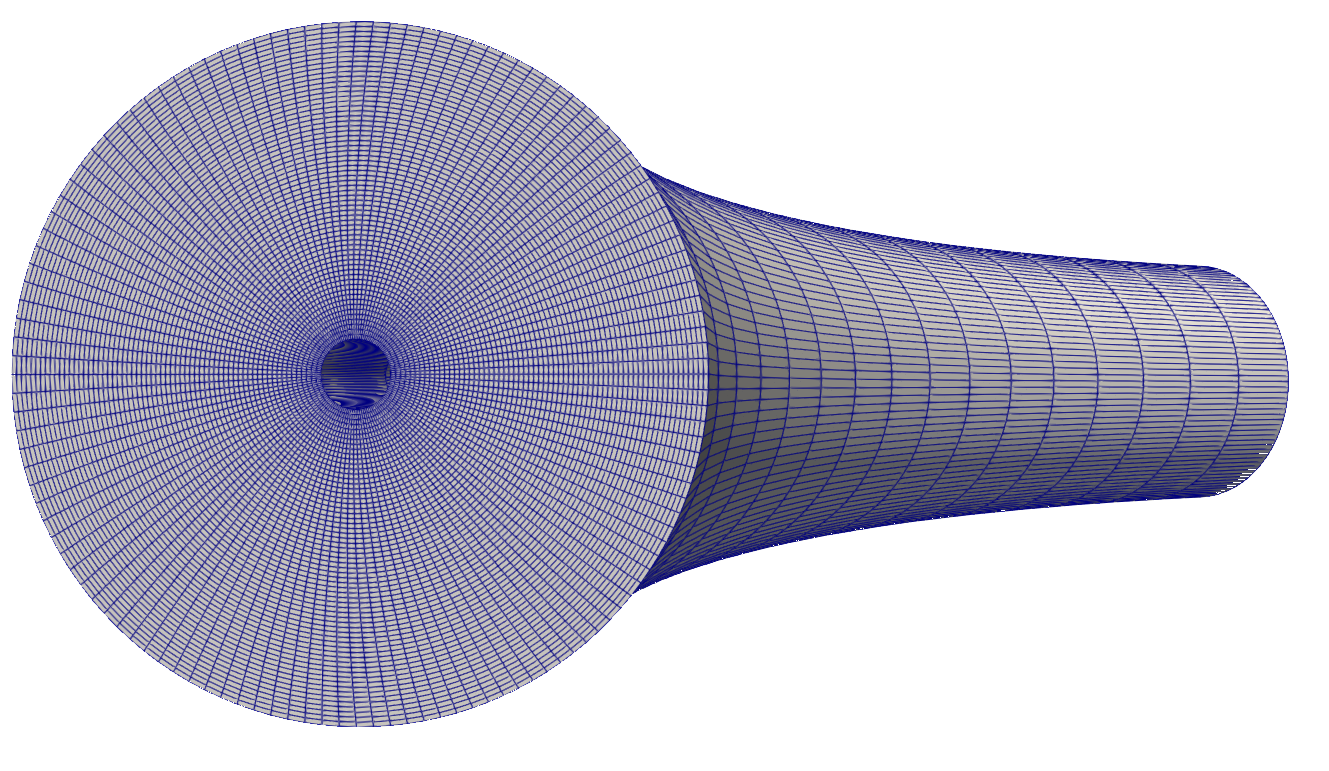}
	\caption{Domain and mesh for the simulation of the TAE}
\label{fig:mesh_tae}
\end{figure}

\begin{figure}
\centering
\captionsetup[subfigure]{justification=centering}
    \begin{subfigure}[b]{0.49\textwidth}
    \centering
    \includegraphics[width=\textwidth]{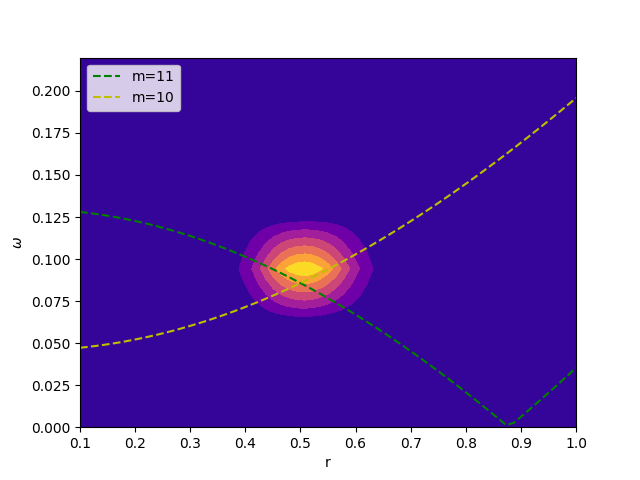}
    \caption{Without artificial viscosity}
    \end{subfigure}
    \begin{subfigure}[b]{0.49\textwidth}
    \centering
    \includegraphics[width=\textwidth]{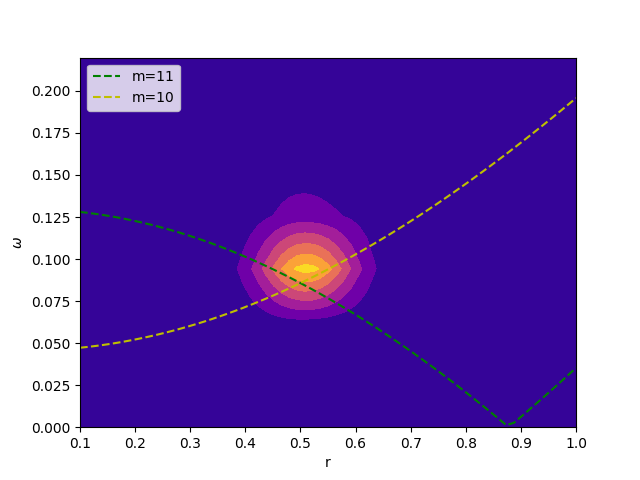}
    \caption{With artificial viscosity}
    \end{subfigure}
\caption{Position of the TAE with respect to the shear alfven frenquency for $m=10$ and $m=11$}
\label{fig:tae_spectrum}
\end{figure}

\begin{figure}
\centering
\captionsetup[subfigure]{justification=centering}
    \begin{subfigure}[b]{0.49\textwidth}
    \centering
    \includegraphics[width=\textwidth]{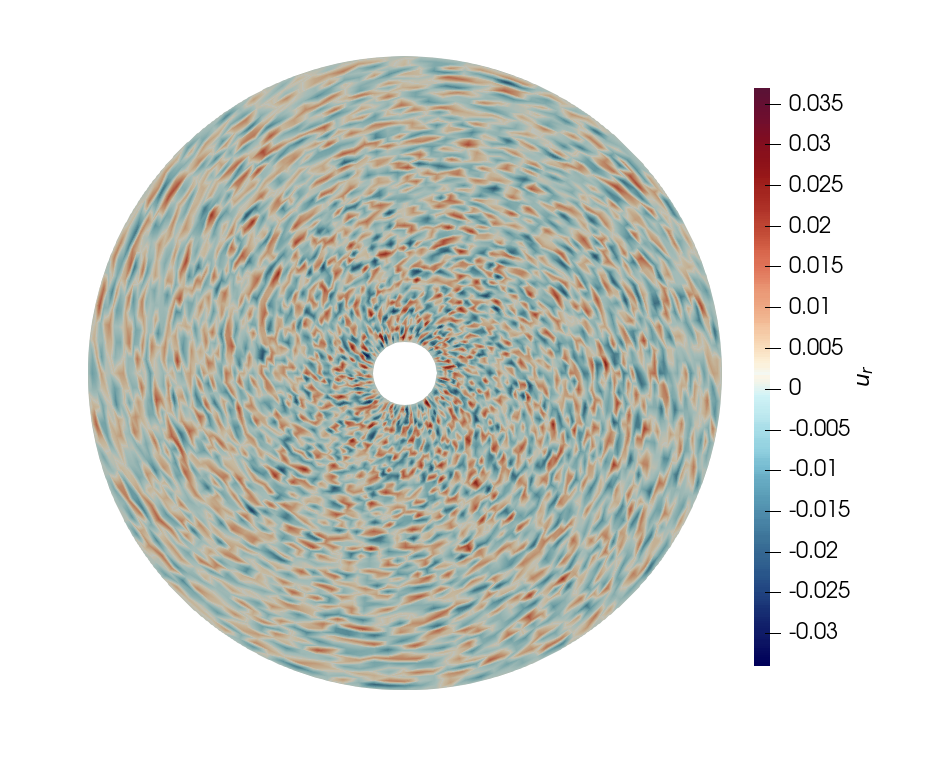}
    \caption{Without artificial viscosity}
    \end{subfigure}
    \begin{subfigure}[b]{0.49\textwidth}
    \centering
    \includegraphics[width=\textwidth]{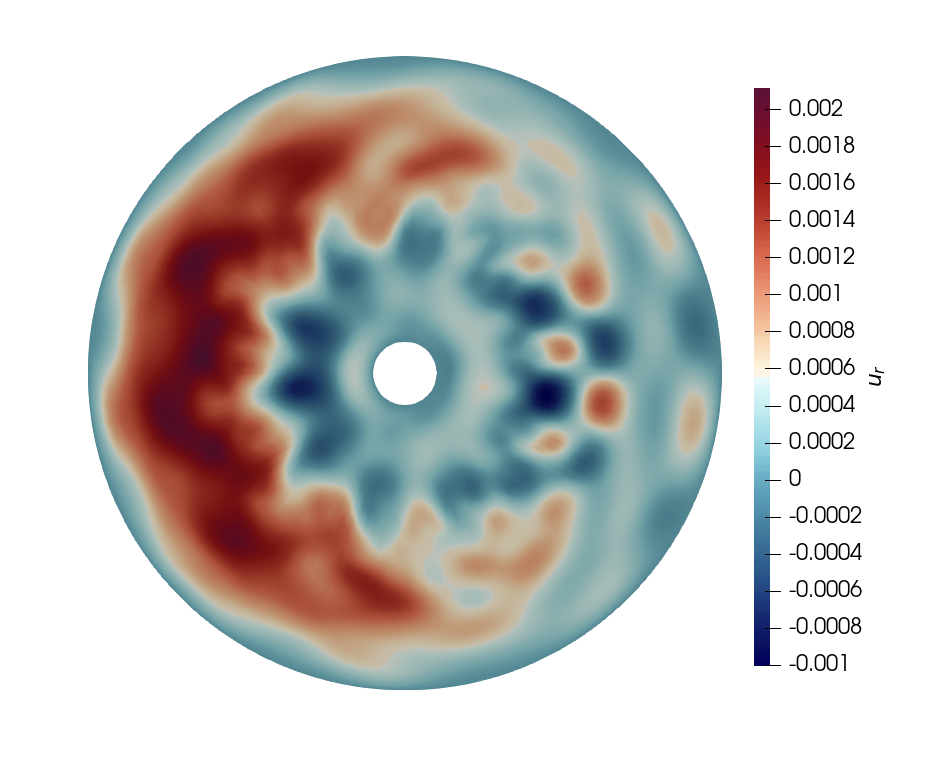}
    \caption{With artificial viscosity}
    \end{subfigure}
\caption{Radial velocity of the TAE simulation at final time, cut at $\phi=0$}
\label{fig:tae_cut}
\end{figure}

Our final ideal test is a Toroidal Alfven Eigenwave. 
This is a mode that can exist in some particular Tokamak configuration, where the frequency of toroidal modes depend on the radial coordinate.  
In this scenario, it can happen that two modes have the same frequency at a certain given radius, creating a particular mode, called Toroidal Alfven Eigenwave. 
We refer to \cite{vlad1999dynamics} for more details on TAEs and plasma waves in Tokamak geometries.
This tests takes place in a simplified Tokamak configuration, a Hollow torus, of major radius $R_0=10.$, and minor radius $r_0=1.$. 
For tests in simplified Tokamak geometries we will use the $(r, \theta, \phi)$ coordinates, where $r$ is the inner radius, $\theta$ the poloidal angle and $\phi$ the toroidal angle.
The hole inside the torus (needed to avoid having to deal with polar singularity, although this could be done using the framework provided by \cite{holderied2021mhd}) has radius of $0.1$. 
The initial conditions are given by the approximate equilibrium:
\begin{align*}
\rho &= 1. ~, \\
p(r) &= \beta  \frac{B_0^2}{2}(p_0-p_1 r^2-p_2 r^4)~, \\
\BB &= \nabla \psi \times \nabla \phi + g \nabla \phi~, \\
\text{with } g &= -B_0 R_0 ~, \\
\text{and } \frac{d \psi}{dr} &= \frac{B_0 r}{q(r) \sqrt{1-r^2/R_0^2}} ~,\\
\text{with } q(r) &= q_0 + (q_1-q_0) r^2
\end{align*}
With parameters: $B_0 = 3.$, $q_0=1.71$, $q_1=1.87$, $p_0=3.$, $p_1 = 0.95$, $p_2=0.05$ and $\beta = 0.002$. 
This initial (approximate) equilibrium is perturbed by $\delta u_r = \epsilon \chi(r)(sin(10 \theta - 6 \phi) + sin(11 \theta - 6 \phi)$ with $\chi(r) = \exp(-(r-0.5)^2/0.01)$ and $\epsilon = 0.01$ which will excite a TAE corresponding to the $10$-th and $11$-th mode in the poloidal direction and $-6$-th mode in the toroidal direction, located around $r=0.5$. 
Due to the symmetry of the problem in the toroidal direction, we only simulate a sixth of the whole torus. 
Simulations are run using $64$ elements in the radial direction, $128$ in the poloidal direction and $16$ in the toroidal direction. We also use field align coordinate, leading to a more bended mesh. We plotted the domain and the mesh for the simulation in \cref{fig:mesh_tae}.
We present two results for this test, with and without artificial stabilization. Both were run with $\Delta t = 0.1$, until a final time of $200$. 
For the stabilized one, we use artificial stabilization parameters $\mu_a = \eta_a = 4.8\times 10^{-5}$, corresponding roughly to $2h^2$ with $h$ the smallest element size. 
In \cref{fig:tae_spectrum} we present a Fourier transform in time of the $\theta$ component of the velocity field, integrated over $\theta$ and cut at $\phi = 0$. 
In this plots we can clearly see the TAE resonance at the frequency where the two branches cross at $r=0.5$. We can also observe that with artificial dissipation the TAE is less well located, but still well captured. 
\cref{fig:tae_cut} presents the radial velocity for cuts of the simulation. We can clearly see on the left figure that without artificial dissipation the small scale oscillations are predominant and make every larger scale structure invisible. On the right plot, we can clearly see the structure of the TAE in the radial velocity component. We can also see that an other wave is superposed (seen by the fact that velocity is positive of the right and negative on the left). This other wave is due to the fact that we only start from an approximate equilibrium and this bigger mode is oscillating around a true equilibrium. However due to the small amplitude of all the oscillations (either small scales one in the non stabilized simulation or bigger ones in the stabilized one), we are still approximately in the linear regime and we can still observe the TAE as shown by \cref{fig:tae_spectrum}.
\subsection{Viscous and Resistive tests}
\label{sec:vr_test}
\subsubsection{One dimensional current sheet}
Our first test is a one dimensional test that has an approximate analytical solution, it is in fact a solution of the vectorial heat equation 
\begin{equation}
\partial_t \BB + \nabla \times (\eta \nabla \times \BB) = 0 ~,
\end{equation}
that is given by 
\begin{equation}
\label{eqn:By_shear}
B_y(x,t) = - B_y^0 \text{erf}(\frac{1}{2} \frac{x}{\sqrt{\eta (t+t_0)}}) ~.
\end{equation}
If we set the background velocity to zero as well as the density and entropy to be constant, it is then a perturbative solution to the VRMHD equations. That means that 
\begin{equation}
\label{eqn:approx_shear}
(\rho, \uu, s, \BB) = (\rho_0, (0,0,0), s_0, (0, B_y(x,t), B_z^0)) ~,
\end{equation}
is an approximate solution to the VRMHD equations. We therefore use the previously mentioned solution to initialize our discrete scheme and let it evolve until $T=1000$ to compare with the approximate given by \cref{eqn:By_shear,eqn:approx_shear}. 
In our numerical experiments, we use $\gamma = 5/3$ (monoatomic perfect gas) and initialize the solution with $\rho_0 = 1.$, $s_0 = 9.62$ (correspond to a constant pressure background with $p_0 = 10^5)$, $B_z^0 = 10^4$ and $B_y^0 = 10^{-3}$. The simulation is run in a one dimensional domain $[-50,50]$, with a time step $\Delta t = 2.\times 10^{-3}$ and $t_0 = 10.$. The magnetic resistivity is set to $\eta = 0.1$. We specificaly use high values of $p_0$ and $B_z^0$ compared to $B_y^0$ to ensure the solution remains in the pertubative domain.

\begin{figure}
    \centering
    \includegraphics[width=0.5\textwidth]{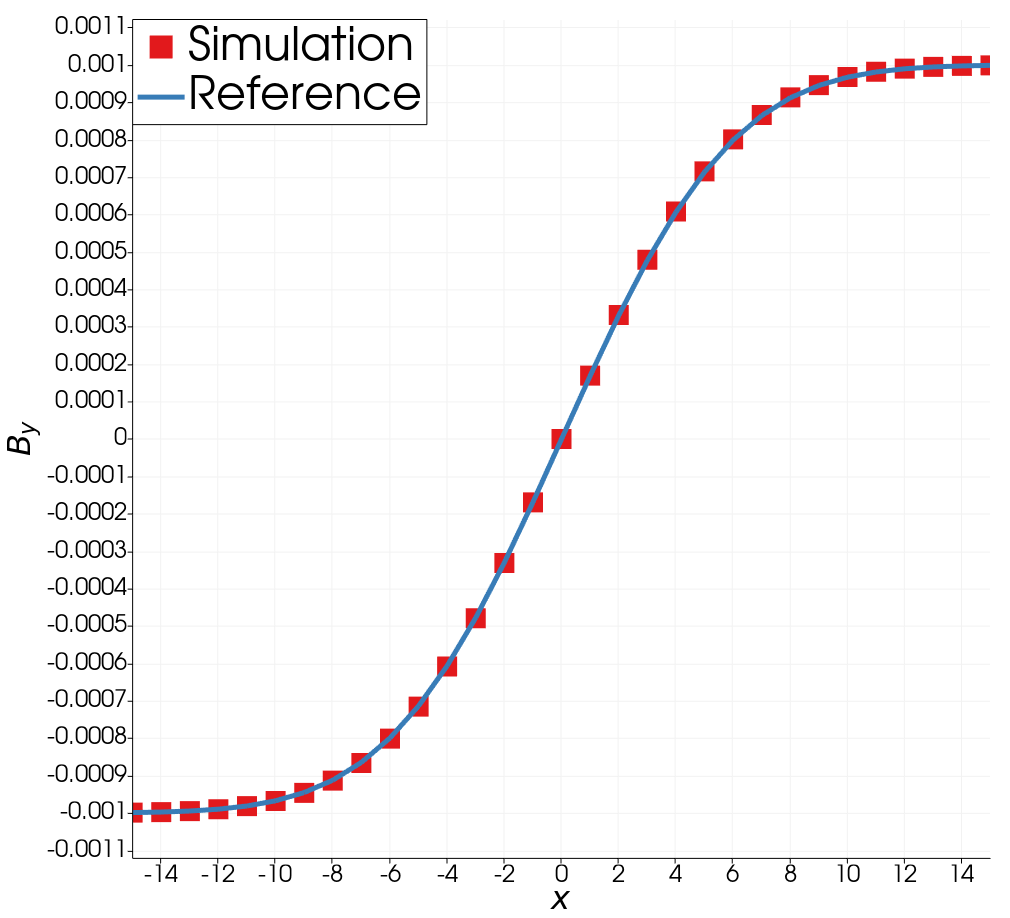}
    \caption{Comparison of the simulated $y$-component of the magnetic field and the reference perturbative solution for the one dimensional current sheet. Zoom on the transition zone around $x=0$.}
    \label{fig:1d_current}
\end{figure}

The results are presented in \cref{fig:1d_current}. We can see almost perfect agreement between the simulated magnetic field and the perturbative solution. This test shows that our scheme is well able to simulate the resistive term.
\subsubsection{Two dimensional current sheet}
Our next test aims at evaluate the ability of our scheme to reproduce plasma instabilities. In this goal we study the evolution of different modes in a perturbed current sheet, trying to reproduce the results obtained in \cite{landi2008three}. 
The simulation take place in a rectangle $[0, 6 \pi] \times [-\frac{\pi}{2}, \frac{\pi}{2}]$, and we use periodic boundary conditions in the $x$ direction. 
The simulated plasma is initialized with the following conditions :
\begin{align*}
\rho(x,y,0) &= 1 ~, \\
s(x,y,0) &= \rho(x,y,0)\log(\frac{5}{2(\gamma-1)}) ~, \\
u_x (x,y,0) &= 0 ~, \\
u_y (x,y,0) &= \epsilon \chi(y) \sum\limits_{k_x} sin(k_x x + \phi_{k_x}) ~,\\
B_x (x,y,0) &= \tanh(y/\delta) ~, \\
B_z (x,y,0) &= 0 ~, \\
B_z (x,y,0) &= \sqrt{1-B_x(x,y,0)^2} ~.
\end{align*} 

Where $\chi(y) = \frac{\tanh(\delta y)}{\cosh(\delta y)}$ and the parameters are $\delta = 0.1$ and $\epsilon = 10^{-4}$, while in the sum, the $k_x$ are the possible wave length $k_x = n/3$ with $n$ positive integer.
If $\epsilon=0$ and there is no resistivity ($\eta = 0$) this setup is a force free current sheet equilibrium. 
We here use $\eta = 2 \times 10^{-4}$, and in this setup, a tearing instability growth in the current sheet. 
The aim here is to compare the growth rate of the different modes, with the ones computed in \cite{landi2008three}. 
However to do so we need to linearize the equation around this equilibrium so that the current sheet is not dissipated by the resistivity. 
This is simply done by replacing the propagator $\Phi^{res}_{\Delta t}$ with 
\begin{subequations}
\begin{equation}
\frac{\Bh^{n+1}- \Bh^n}{\Delta t} + \nabla \times(\eta \tilde{\nabla} \times (\Bh^{n+1} - \Bh^0)) = 0  ~,
\end{equation} 
\begin{equation}
\int_\Omega \frac{\rhoh^n e(\rhoh^n,\sh^{n+1})-\rhoh^n e(\rhoh^n,\sh^n)}{\Delta t} q_h - \int_\Omega \eta \nabla \times \Bh^{n+\frac{1}{2}} \cdot \nabla \times (\Bh^{n+1} - \Bh^0) q_h = 0 \qquad \forall q_h \in V^3_h ~.
\end{equation}
\end{subequations}
\begin{figure}
\centering
\captionsetup[subfigure]{justification=centering}
    \begin{subfigure}[b]{0.49\textwidth}
    \centering
    \includegraphics[width=\textwidth]{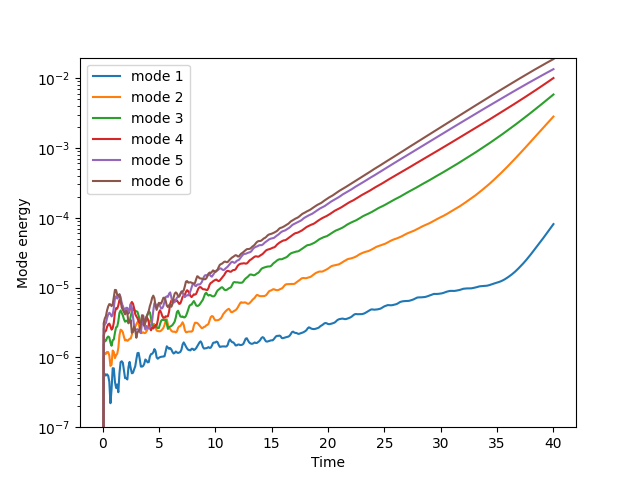}
    \caption{Growth of the energy carried by the 6 first modes.}
    \end{subfigure}
    \begin{subfigure}[b]{0.49\textwidth}
    \centering
    \includegraphics[width=\textwidth]{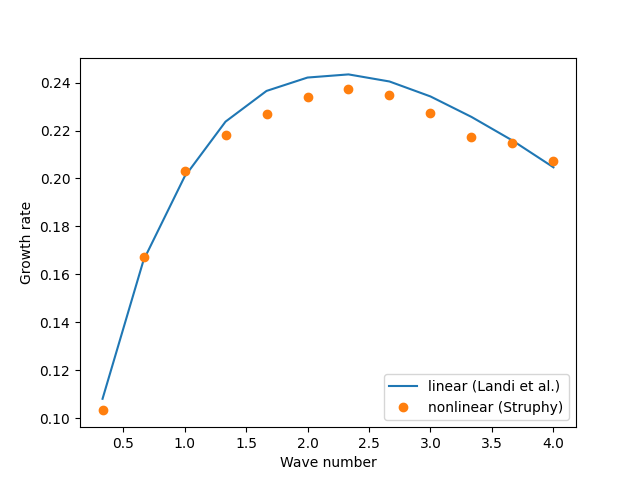}
    \caption{Growth rate obtained by struphy and comparison with reference.}
    \end{subfigure}
\caption{Results of the tearing instability}
\label{fig:current_sheet}
\end{figure}
Results are shown in \cref{fig:current_sheet}, those are obtained by doing a Fourier decomposition of the magnetic field and compute the energy in each of the modes. 
For this simulation we used a $128 \times 256$ grid, a maximal spline degree of $p=2$ and constant time step $\Delta t = 0.1$.
We can clearly observe that in a first phase (time $0-15$) some noise in produced, probably because of numerical instability, then in a second phase (time $15-30$), we have linear growth of every mode and around time $35$ the nonlinearity start to couple the modes, ending the linear growth phase. 
The growth rate plotted correspond to the ones computed during the second phase. 
We can see that those are in great agreement with the ones from the reference even if we use a non-linear scheme, while they were computed using a linear approximation. 
This result gives us some great confidence that our variational scheme is able to reproduce well plasma instabilities.
\subsubsection{Viscoresistive Orszag-Tang Vortex}
We now study the evolution of a viscous and resistive Orszag-Tang vortex \cite{warburton1999discontinuous}. As in the ideal case, the simulation take place in a periodic box $[0, 2\pi]^2$, and the initial condition are given by:
\begin{align*}
\rho(x,y,0) &= 1. ~, \\
s(x,y,0) &= \log\Big(\frac{p}{\gamma-1}\Big) ~, \\
p(x,y,0) &= \frac{15}{4} + \frac{1}{4} \cos(4x) +  \frac{4}{5}\cos(2x)\cos(y)-\cos(x)\cos(y)+\frac{1}{4}\cos(2y) ~, \\
\uu (x,y,0) &= (-\sin(y), \sin(x)) ~,\\
\BB (x,y,0) &= (-\sin(y), \sin(2x)) ~ 
\end{align*}
For this test we use $\gamma = \frac{5}{3}$ and the dissipative parameters are $\eta = \mu = 0.01$
\begin{figure}
\centering
\captionsetup[subfigure]{justification=centering}
    \begin{subfigure}[b]{0.49\textwidth}
    \centering
    \includegraphics[width=\textwidth]{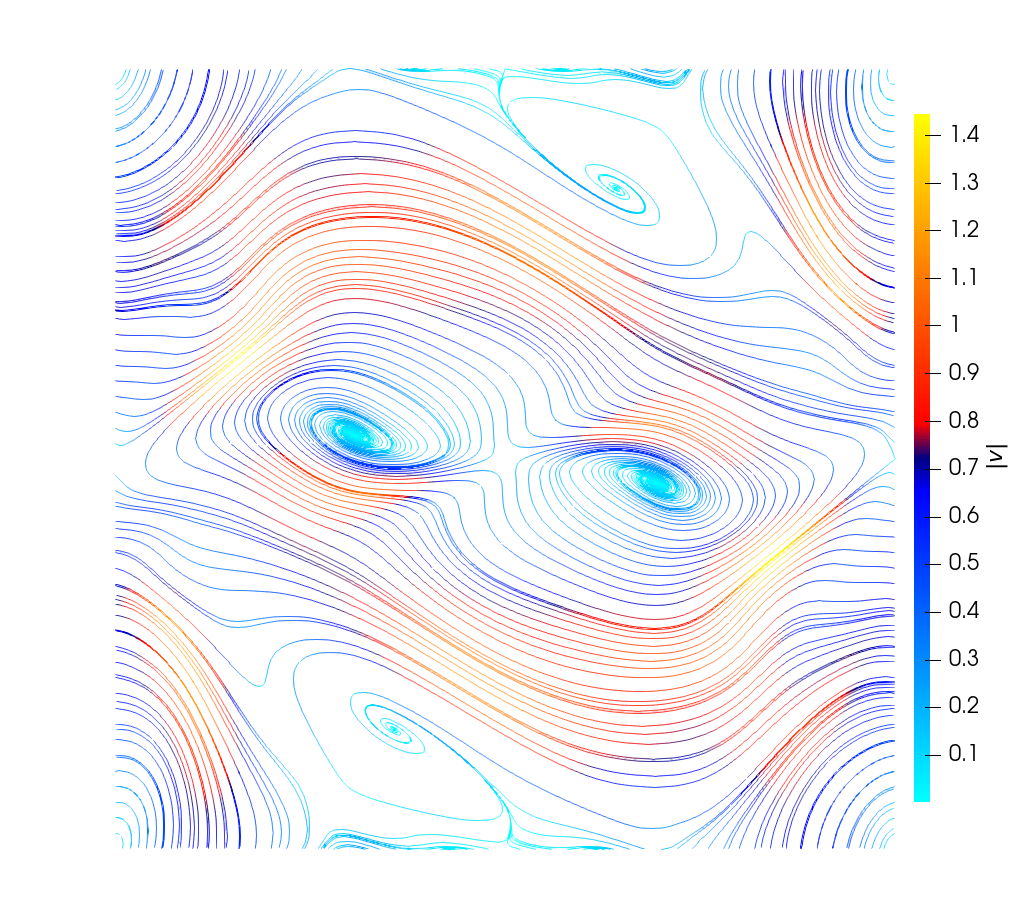}
    \caption{Streamlines of the velocity coloured by the $L^2$ norm of the velocity}
    \end{subfigure}
    \begin{subfigure}[b]{0.49\textwidth}
    \centering
    \includegraphics[width=\textwidth]{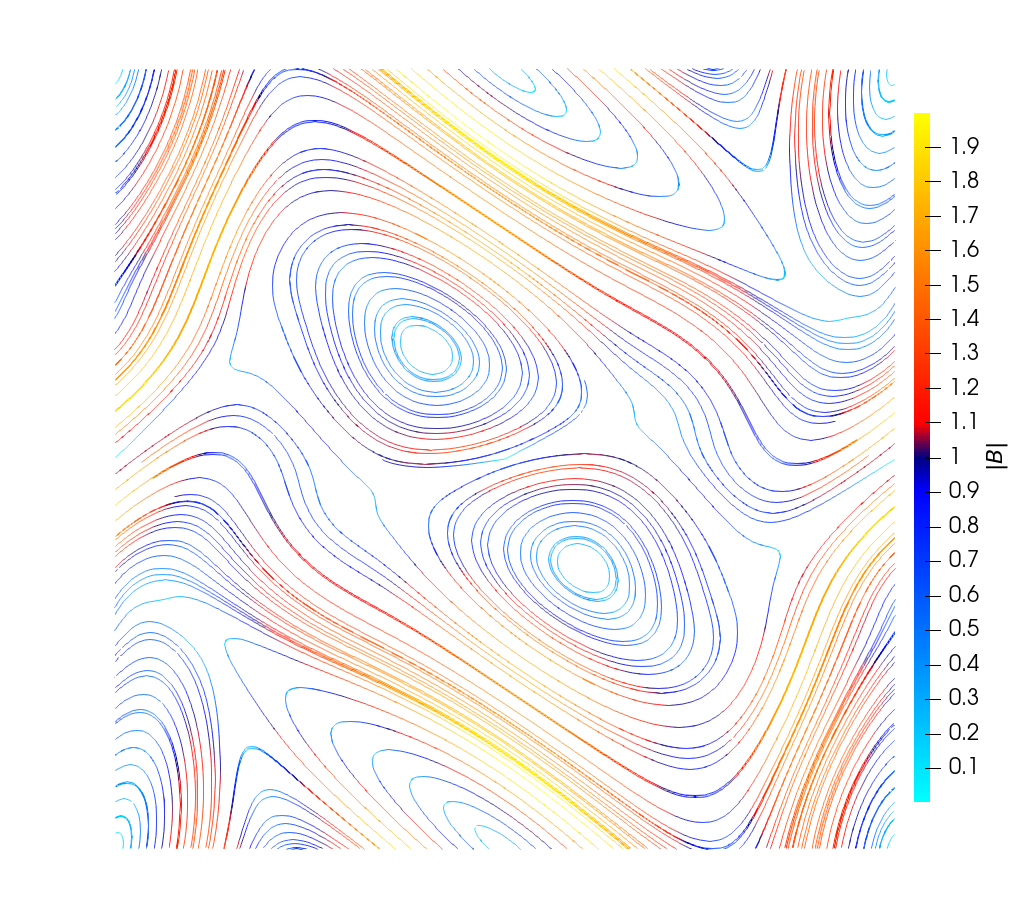}
    \caption{Streamlines of the magnetic field coloured by the $L^2$ norm of the magnetic field}
    \end{subfigure}
\caption{Viscoresistive Orszag-Tang vortex at $t=2$}
\label{fig:VROT}
\end{figure}
We present the results for a $256 \times 256$ grid, using splines of maximal degree $p=2$ in \cref{fig:VROT}, using a constant time step $\delta t = 10^{-3}$. We plot the streamlines of the $\uu$ and $\BB$ fields, at time $t=2.$, and can observe that our results are in great accordance with the literature \cite{warburton1999discontinuous,dumbser2009high,fambri2021novel}. This test allows us to check that all the terms, specially the viscous and resistive one are well integrated by our variational scheme, specially that the field lines of $\BB_h$ remain closed even with the addition of resistivity.
\subsubsection{Resistive kink mode}
Our last test is a resisitve kink mode, which is a unstable mode existing for certain plasma equilibrium when resistivity is included in the MHD model.
This test in ran in a simplified tokamak geometry, similar to the TAE test \cref{sec:tae_test} but with the full toroidal domain being discretized (not only one sixth as before) and using field align coordinates (not straight line going from the pole to the edge but coordinates following the magnetic field). For this test we needed to substract the equilibrium flow to the equation, in order to balance the fact that our starting equilibrium is not a discrete equilibrium of our scheme. This is simply done by solving $\partial_t U + F_h(U) = F_h(U_0)$ where $F_h$ denotes our scheme, $U$ a generic unknown (here the vector of all our unknown) and $U_0$ is our starting approximate equilibrium. Although this breaks the structure preservation of our approach, the source terms are relatively small and we are still able to have a very good preservation of the key properties. Future research would focus on tackling this problem.
The plasma is initialized with an equilibrium provided by the GVEC MHD equilibrium solver \cite{hindenlang_gvec_2019}, with on axis magnetic field $\BB_0 = 1.$, density $\rho_0 = 1.$ and pressure $p_0 = 2e-3$ and on the edge $\rho_e = 0.1 \rho_0$ and $p_e = 0.02 p_0$, similar to the equilibrium used in \cite{haverkort2016implementation}. This equilibrium is perturbed by random white noise added in order to excite all frequencies, with norm $1e-8$. For the dissipative parameters, we use a relatively high resistivity of $\eta = 1e-4$ in order to have a relatively large kink mode (the size of the mode and the growth rate are proportional to $\eta^{1/3}$) and viscosity of $\mu = 1e-6$. The discretization uses $32 \times 16 \times 16$ elements, a maximal spline degree of $p=2$ and a relatively large time step of $\Delta t = 1e-1$. 

\begin{figure}
\centering
\captionsetup[subfigure]{justification=centering}
    \begin{subfigure}[b]{0.49\textwidth}
    \centering
    \includegraphics[width=\textwidth]{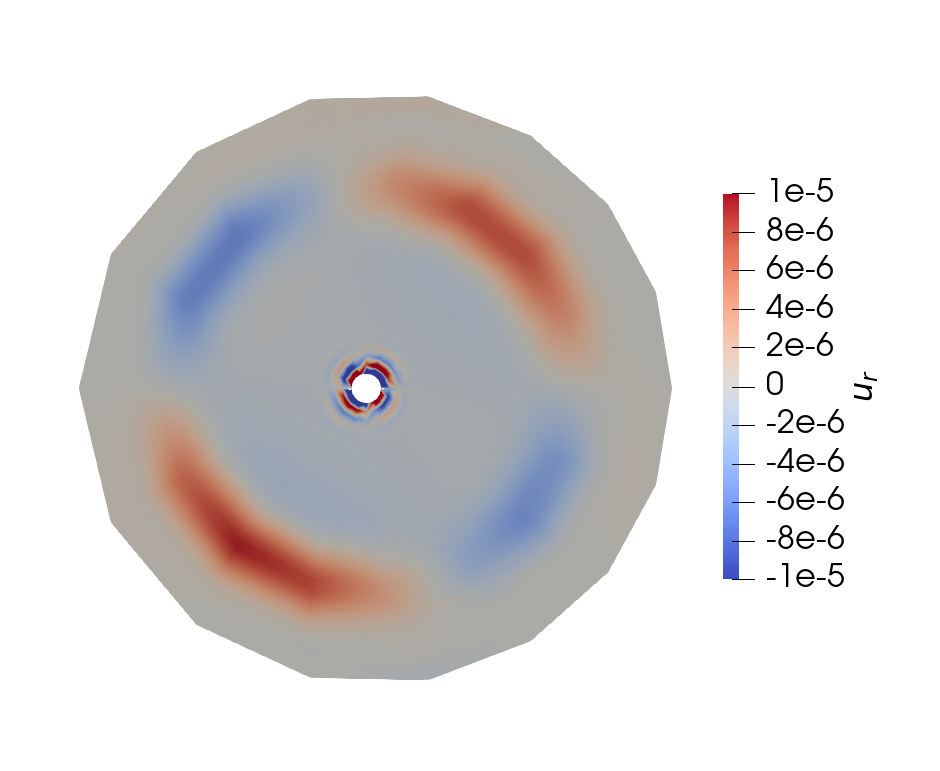}
    \caption{$R$ component of the velocity field}
    \end{subfigure}
    \begin{subfigure}[b]{0.49\textwidth}
    \centering
    \includegraphics[width=\textwidth]{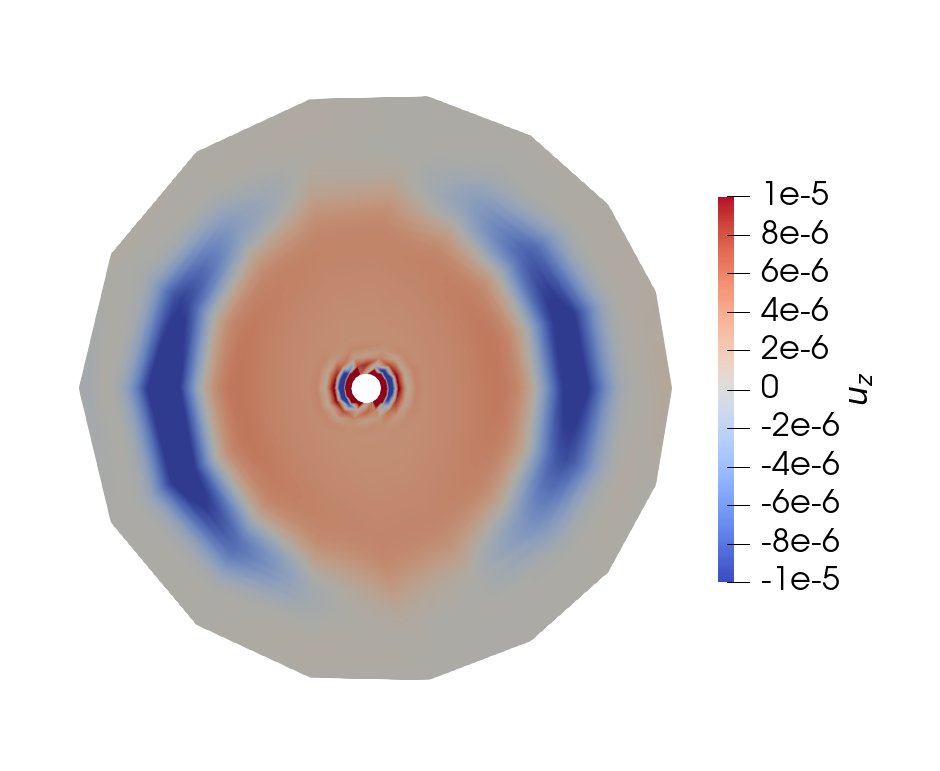}
    \caption{$Z$ component of the velocity field}
    \end{subfigure}
\caption{Kink mode at $t=2000.$}
\label{fig:Kink_u}
\end{figure}

\begin{figure}
\centering
    \includegraphics[width=\textwidth]{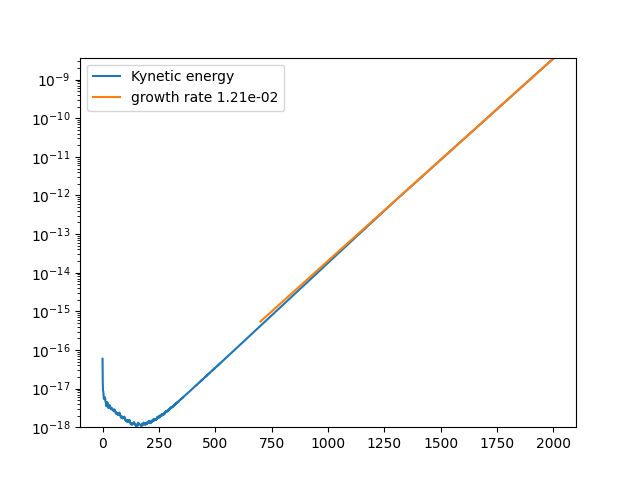}
    \caption{Kynetic energy growth of the kink mode}
\label{fig:Kink_growth}
\end{figure}

\Cref{fig:Kink_u} presents the velocity fields at $t=2000.$ for a toroidal cut. We can clearly recognize the characteristic form of the kink mode \cite{haverkort2016implementation}, with on this cut the central velocity going upward. In \cref{fig:Kink_growth} we present the growth of the kynetic energy, which is due to the kink mode. We clearly see a linear growth, showing that even if our code is non-linear and well able to resolve non-linear dynamics, we are still able to have a very good accordance for the linear phase of instability growth. The growth rate for the kink mode is here $1.21e-2$ which is in the range of the results found in the literature \cite{kerner1998castor}. This tests clearly show the ability of our approach to well reproduce linear instabilities, using bigger time steps then usual approach (we here use $\Delta t = 1e-1$). However we also see some spurious features around the axis due to the use of a hollow torus geometry. Future work will try to overcome this limit by using the polar-spline framework, which was unfortunately not easily adaptable for this work due to the transformation of the velocity field, not fitting into previous approaches. 

\section{Conclusion and perspective}
\label{sec:concl}

In this work we presented a new scheme for solving the equation of viscous and resistive magnetohydrodynamics. Our approach is based on a discrete variational principle mimicking the continuous one, and allow for a algorithm that naturally preserves invariants from the system. We were able to show that the obtained semi-discrete scheme also has a metriplectic structure, and we propose a time splitting algorithm for the time integration that enjoy this structure. 

We then implemented our approach in the Struphy library and conduced intensive testing, where we used our framework first to stabilize ideal simulations, greatly improving the results from the non-stabilized version. Second we tested our algorithm on plasma instability triggered by resistivity and were able to obtain satisfying results. Overall the results show that variational schemes could compete with more traditional approach for the study of instability, although a lot of improvement should still be done.

Future work on this topic will focus on solving the limitation pointed in the numerical section: simulation of the magnetic axis (using polar splines) and simulation around equilibrium (to overpass the need for subtracting the equilibrium flow). We would also like to improve the scalability of this approach in order to compete with established MHD codes. Some interested is also put in developing a delta-f variational approach that could lead to better result while being less computationally intensive. 

\bibliographystyle{plain} 
\bibliography{VRMHD}

\end{document}